\documentclass{amsart}
\usepackage{a4wide}
\usepackage{amssymb,amsmath,amscd, amsfonts,latexsym}
\usepackage{graphicx}
\usepackage{amssymb}
\usepackage{epstopdf}
\usepackage{verbatim}  

\usepackage[all]{xy}
\usepackage{color}

\newtheorem{theorem}{Theorem}

\newtheorem{lemma}{Lemma}
\newtheorem{proposition}{Proposition}
\theoremstyle{remark}
\newtheorem{remark}{Remark}

\newcommand{\bbQ}{\mathbb Q}

\newcommand{\bbZ}{\mathbb Z}
\newcommand{\bbN}{\mathbb N}

\title[Sums of Fibonacci numbers and powers of $2$]{On Diophantine equations involving sums of Fibonacci numbers and powers of $2$} 

\author[K. C. Chim]{Kwok Chi Chim}

\author[V. Ziegler]{Volker Ziegler}

\thanks{The first author was supported by the Austrian Science Fund (FWF) under the projects P26114 and W1230.}



\address{K. C. Chim \newline
         \indent Institute of Analysis and Number Theory, Graz University of Technology \newline
        \indent Kopernikusgasse 24/II \newline
         \indent A-8010 Graz, Austria}
        
\email{chim\char'100math.tugraz.at}

\address{V. Ziegler \newline
         \indent University of Salzburg \newline
         \indent Hellbrunnerstrasse 34/I \newline
         \indent A-5020 Salzburg, Austria}

\email{volker.ziegler\char'100sbg.ac.at}

\allowdisplaybreaks[1]

\begin{document}

\begin{abstract}
In this paper, we completely solve the Diophantine equations $F_{n_1} + F_{n_2} = 2^{a_1} + 2^{a_2} + 2^{a_3}$ and $  F_{m_1} + F_{m_2} + F_{m_3} =2^{t_1} + 2^{t_2} $,
where $F_k$ denotes the $k$-th Fibonacci number. In particular, we prove that $\max \{n_1, n_2, a_1, a_2, a_3 \}\leq 18$ and $\max \{ m_1, m_2, m_3, t_1, t_2 \}\leq 16$.
\end{abstract}

\maketitle

\section{Introduction}\label{Sec:Intro}

There is a vast literature on solving Diophantine equations involving the sequence $\{F_n\}_{n\geq 0}$ of Fibonacci numbers (defined by $F_0=0$, $F_1=1$ and $F_{n+2}=F_{n+1}+F_n$ for $n\geq 0$),
the sequence $\{F_n^{(k)}\}_{n\geq 0}$ of $k$-generalized Fibonacci numbers, the sequence $\{P_n\}_{n\geq 0}$ of Pell numbers 
or other recurrence sequences. For instance, recent results include Bravo and Luca \cite{Luca16} where they studied 
the Diophantine equation 
$$F_n+F_m = 2^a.$$
In \cite{Luca14b} they extended their work to $k$-generalized Fibonacci number $F_n^{(k)}$, and studied 
the equation 
$$F_n^{(k)} + F_m^{(k)} = 2^a.$$
Besides, Bravo, Faye and Luca \cite{BBL16} studied the Diophantine equation 
$$P_l + P_m + P_n = 2^a.$$

The most general results in this respect are due to Stewart \cite{Stewart:1980}, who studied representations of integers in two different bases. Note that e.g. the result due to Bravo and Luca \cite{Luca16}
can be seen as an attempt to find all integers that have only few digits in base $2$ and the Zeckendorf expansion simultaneously. Also Luca \cite{Luca:2000} proves a similar result.
Finally let us mention a recent result due to Meher and Rout \cite{Meher:2016} on the Diophantine equation
$$U_{n_1}+\dots + U_{n_t}=b_1 p_1^{z_1}+\dots +b_s p_s^{z_s}$$
in non-negative integers $n_1 ,\dots , n_t , z_1 , \dots , z_s$, where $\{U_n\}_{n\geq 0}$ is a binary, non-degenerate recurrence sequence with positive discriminant, $b_1,\dots,b_s$
are fixed non-negative integers and $p_1,\dots,p_s$ are fixed primes.

Also recently Diophantine equations have been studied which can be regarded as variants of Pillai's problem \cite{Pillai:1936}. For instance, 
Chim, Pink and Ziegler \cite{SalzburgI} obtained all the integers $c$ such that the Diophantine equation
$$F_n - T_m=c$$
has at least two solutions. Here $T_m$ denotes the $m$-th Tribonacci number. Ddamulira, Luca, and Rakotomalala \cite{Luca15} considered 
the Diophantine equation
$$F_n - 2^m=c$$
and found all integers $c$ for which this Diophantine equation has at least two solutions.
Recently, Bravo, Luca and Yaz\'an \cite{Bravo16} considered the Diophantine equation
$$T_n - 2^m=c$$
instead. The most general result is due to
Chim, Pink and Ziegler \cite{SalzburgII} who considered the case, where $U_n$ and $V_m$ are the $n$-th and $m$-th numbers in linear recurrence sequences
$\lbrace U_n\rbrace_{n\geq 0}$ and $\lbrace V_m\rbrace_{m\geq 0}$ respectively and found effective upper bounds for $|c|$ such that the Diophantine equation
$$U_n-V_m=c$$
has at least two solutions.

All the problems stated above are solved by a similar strategy, the iterated application of linear forms in logarithms. We extend this strategy and study the 
two Diophantine equations 
$$F_{n_1} + F_{n_2} = 2^{a_1} + 2^{a_2} + 2^{a_3}$$
and
$$F_{m_1} + F_{m_2} + F_{m_3} =2^{t_1} + 2^{t_2}.$$
In particular, we prove the following two theorems.

\begin{theorem}\label{Main3A_1}
Let $(n_1,n_2,a_1,a_2,a_3)\in \bbN^5$ be a solution to the Diophantine equation 
\begin{equation}  \label{eqt3A_1}
F_{n_1} + F_{n_2} = 2^{a_1} + 2^{a_2} + 2^{a_3}
\end{equation} 
such that $n_1 \geq n_2 \geq 0$ and $a_1 \geq a_2 \geq a_3 \geq 0$,  then $n_1 \leq 18$ and $a_1\leq 11$. In particular,
equation \eqref{eqt3A_1} has exactly $78$ solutions.
\end{theorem}

\begin{theorem}\label{Main3A_2}
Let $(m_1, m_2, m_3, t_1, t_2) \in \bbN^5$ be a solution to the Diophantine equation 
\begin{equation}  \label{eqt3A_2}
F_{m_1} + F_{m_2} + F_{m_3} = 2^{t_1} + 2^{t_2} 
\end{equation} 
such that $m_1 \geq m_2  \geq m_3 \geq 0$ and $t_1 \geq t_2 \geq 0$, then $m_1 \leq 16$ and $t_1\leq 10$. In particular,
equation \eqref{eqt3A_2} has exactly $116$ solutions.
\end{theorem}

\begin{remark}
 The list of solutions to equations \eqref{eqt3A_1} and \eqref{eqt3A_2} is given in the Appendix.
 So we keep the statement of Theorems \ref{eqt3A_1} and \ref{eqt3A_2} short and compact. 
\end{remark}

We shall prove both Theorems \ref{Main3A_1} and \ref{Main3A_2} by the typical strategy also performed in \cite{Luca16,Bravo16,SalzburgI,SalzburgII,Luca15}.
First we extract by a simple computer search all solutions $(n_1, n_2, a_1, a_2, a_3)$ with $n_1 < 360$ to equation \eqref{eqt3A_1} and all solutions $(m_1, m_2, m_3, t_1, t_2)$ with $m_1 < 360$
to equation \eqref{eqt3A_2}, respectively. The key argument to obtain upper bounds for $n_1=\max \lbrace n_1, n_2, a_1, a_2, a_3 \rbrace$ and $m_1=\max \lbrace m_1, m_2, m_3, t_1, t_2 \rbrace$
respectively is to apply lower bounds for linear forms in logarithms. This is done in the seven steps described below, where $c_1,\dots,c_7$ denote effectively computable constants. These
seven steps are in case of the proof of Theorem \ref{Main3A_1} the following:

\begin{description}
 \item[Step 1] We obtain an upper bound 
 $$\min\{(a_1-a_2)\log 2,(n_1-n_2)\log \alpha\}\leq c_1 \log n_1.$$
 Hence we have to distinguish between the following two cases:
 \begin{description}
  \item[Case 1] $\min\{(a_1-a_2)\log 2,(n_1-n_2)\log \alpha\}=(a_1-a_2)\log 2\leq c_1 \log n_1$ 
  \item[Case 2] $\min\{(a_1-a_2)\log 2,(n_1-n_2)\log \alpha\}=(n_1-n_2)\log \alpha\leq c_1 \log n_1$
 \end{description}
 \item[Step 2] We consider Case 1 and show that $(a_1-a_2)\log 2\leq c_1 \log n_1$ yields 
 $$\min\{(a_1-a_3)\log 2,(n_1-n_2)\log \alpha\}\leq c_2 (\log n_1)^2.$$
 Thus we have to further subdivide Case 1 into the following two cases:
 \begin{description}
  \item[Case 1A] $\min\{(a_1-a_3)\log 2,(n_1-n_2)\log \alpha\}=(a_1-a_3)\log 2\leq c_2 (\log n_1)^2$ 
  \item[Case 1B] $\min\{(a_1-a_3)\log 2,(n_1-n_2)\log \alpha\}=(n_1-n_2)\log \alpha\leq c_2 (\log n_1)^2$
 \end{description}
 \item[Step 3] We consider Case 1A and show that $(a_1-a_3)\log 2\leq c_2 (\log n_1)^2$ implies that
 $$(n_1-n_2)\log \alpha\leq c_3 (\log n_1)^3.$$
 \item[Step 4] We consider Case 1B and show that $(a_1-a_2)\log 2\leq c_1 \log n_1$ and $(n_1-n_2)\log \alpha\leq c_2 (\log n_1)^2$ yield the upper bound 
 $$(a_1-a_3)\log 2\leq c_4(\log n_1)^3.$$
 \item[Step 5] We consider Case 2 and show that $(n_1-n_2)\log \alpha\leq c_1 \log n_1$ yields the upper bound 
 $$(a_1-a_2)\log 2\leq c_5(\log n_1)^2.$$
 \item[Step 6] We continue to consider Case 2 and show that assuming the upper bounds $(a_1-a_2)\log 2\leq c_5(\log n_1)^2$ and $(n_1-n_2)\log \alpha\leq c_1 \log n_1$
 yield the upper bound 
 $$(a_1-a_3)\log 2 \leq c_4(\log n_1)^3.$$
 This is basically Step 4 again, but with probably slightly different constants.
 However after Step 6 we have found upper bounds for $(a_1-a_2)\log 2$, $(a_1-a_3)\log 2$ and $(n_1-n_2)\log \alpha$.
 \item[Step 7] We show that the upper bounds found in the previous steps yield an inequality of the form $n_1\leq c_7 (\log n_1)^4$. Thus we obtain an absolute bound for $n_1$. 
\end{description}

As soon as we have found an upper bound for $n_1$ we go through all seven steps again but apply instead of lower bounds for linear forms in logarithms the Baker-Davenport reduction
method and obtain in all steps small, absolute bounds respectively. In case the Baker-Davenport reduction method fails we can make use of a criteria of Legendre for continued fractions to reduce the huge upper bounds
to rather small upper bounds. Indeed we succeed to show that all solutions satisfy $n_1<360$, which already have been found by our previous computer search.

Of course, a slight modification of these seven steps also leads to a proof of Theorem \ref{Main3A_2}.

It should be noted that due to having more terms in each equation as compared to the equations considered in \cite{Luca16,Bravo16,SalzburgI,SalzburgII,Luca15}, we apply several times more the results 
of linear forms in logarithms and the reduction method. E.g. instead of using only twice the results on linear forms in logarithms and the reduction method as in \cite{Luca16} we apply them seven times.

\section{Preliminaries}

In this section, the result of linear forms in logarithms by Baker and W\"ustholz~\cite{bawu93} is stated.
Besides, we state a lemma from \cite{BBL16}, which is 
a generalization of a result due to Baker and Davenport~\cite{BD69} the so-called Baker-Davenport reduction method. Both results will be used to prove Theorems \ref{Main3A_1} and \ref{Main3A_2}.

\subsection{A lower bound for linear forms in logarithms of algebraic numbers}

In 1993, Baker and W\"ustholz \cite{bawu93} obtained an explicit bound for linear forms in logarithms with a linear dependence on $\log B$, where $B \geq e$ denotes an upper bound for the height of the linear form (to be defined later in this section). 
It is a vast improvement compared with lower bounds with a dependence on higher powers of $\log B$ in preceding
publications by other mathematicians in particular Baker's original results \cite{Baker:1966}.

Denote by $\alpha_1, \dots, \alpha_k$ algebraic numbers, not $0$ or $1$, and by $\log \alpha_1, \dots, \log \alpha_k$
a fixed determination of their logarithms. Let $K=\bbQ(\alpha_1, \dotso, \alpha_k)$ and let $d=[K:\bbQ]$ be the degree of $K$ over $\bbQ$.
For any $\alpha \in K$, suppose that its minimal polynomial over the integers is
\[ g(x) = a_0 x^{\delta} + a_1 x^{\delta - 1} + \cdots + a_{\delta} = a_0 \prod_{j=1}^{\delta} \left(x - \alpha^{(j)}\right)\]
where $\alpha^{(j)},\; j=1, \dotso, \delta$, are all the roots of $g(x)$. The absolute logarithmic Weil height of $\alpha$ is defined as
\[ h_0(\alpha) = \frac{1}{\delta} \left( \log |a_0| + \sum_{j=1}^{\delta} \log \left( \max \left\lbrace|\alpha^{(j)}|, 1 \right\rbrace \right)\right). \]
Then the modified height $h'(\alpha)$ is defined by
\[ h'(\alpha)=\frac{1}{d}\max \{h(\alpha), |\log\alpha|, 1\},\]
where $h(\alpha) = d h_0(\alpha)$ is the standard logarithmic Weil height of $\alpha$.

Let us consider the linear form
\[ L(z_1, \dotso, z_k)=b_1 z_1 + \cdots + b_k z_k, \]
where $b_1, \dotso, b_k$ are rational integers, not all 0 and define
\[ h'(L) = \frac{1}{d}\max\{h(L),1\}, \]
where $h(L) = d \log \left(\max_{1 \le j \le k} \left\{\frac{|b_j|}{b}\right\}\right)$ is the logarithmic Weil height of $L$,
with $b$ as the greatest common divisor of $b_1, \dotso, b_k$.
If we write $B=\max\lbrace|b_1|, \dotso, |b_k|, e \rbrace$, then we get
\[ h'(L)\leq \log B. \]

With these notations we are able to state the following result due to Baker and W\"ust\-holz~\cite{bawu93}.

\begin{theorem}\label{BaWu}
If ${\it\Lambda}=L(\log \alpha_1, \dotso, \log\alpha_k) \neq 0$, then
\[\log|{\it\Lambda}|\geq -C(k,d)h'(\alpha_1)\cdots h'(\alpha_k)h'(L), \]
where
\[ C(k,d)=18(k+1)!\,k^{k+1}(32d)^{k+2}\log(2kd). \]
\end{theorem}

With $|{\it\Lambda}| \leq \frac{1}{2}$, we have $\frac{1}{2}|{\it\Lambda}| \leq |{\it\Phi}| \leq 2|\it{\Lambda}|$, where
\[ {\it\Phi} = e^{{\it\Lambda}}-1 = \alpha_1^{b_1} \cdots \alpha_k^{b_k}-1, \]
so that
$$
\log\left|\alpha_1^{b_1} \cdots \alpha_k^{b_k}-1\right| \geq \log|{\it\Lambda}| - \log2.
$$

We apply Theorem \ref{BaWu} mainly in the situation where $K=\bbQ(\sqrt{5})$, $k=3$ and $d=2$. In this case we obtain
\[C(3,2)=18 \cdot 4!\cdot 3^4 \cdot 64^5 \log 12 < 9.34 \cdot 10^{13}.\]
We will use this value throughout the paper without any further reference. Besides, let us recall some well known properties of the absolute logarithmic height:
\begin{equation*}
\begin{split}
h_0(\eta \pm \gamma) \leq & \ h_0(\eta)+h_0(\gamma)+\log{2},  \\
h_0(\eta\gamma^{\pm 1}) \leq & \ h_0(\eta)+h_0(\gamma),  \\
h_0(\eta^\ell)=& \ |\ell|h_0(\eta),
\end{split}
\end{equation*}
where $\eta,\gamma$ are some algebraic numbers and $\ell \in \bbZ$.


\subsection{A generalized result of Baker and Davenport}
The following result will be used to reduce the huge upper bounds for $n_1$ and $m_1$ found in Propositions \ref{prop:bound3A_1} and \ref{prop:bound3A_2} respectively. 
Let us state Lemma 6 in~\cite{BBL16} which is regarded as  
a generalization of a result due to Baker and Davenport~\cite{BD69}.
We denote by $\|x\| = \min \{ |x-n| : n \in \bbZ \}$ the distance from $x\in \mathbb R$ to the nearest integer.

\begin{lemma} \label{Reduction}
Let $M$ be a positive integer, let $p/q$ be a convergent of the continued fraction of the irrational $\gamma$ such that $q > 6M$,
and let $A, B, \mu$ be some real numbers with $A > 0$ and $B > 1$. Let $\varepsilon := \|\mu q\| - M \| \gamma q\|$.
If $\varepsilon > 0$, then there is no solution to the inequality
\begin{equation}\label{eq:Baker-Davenport} 0 < \left| u \gamma - v + \mu \right| < AB^{-w}, \end{equation}
in positive integers $u, v$ and $w$ with
\[ u \leq M \quad \quad and \quad \quad w \geq \frac{\log (Aq/ \varepsilon )}{\log B}. \]
\end{lemma}

\begin{remark} \label{RemarkReduction}
Let us explain how we will make use of Lemma \ref{Reduction} and explain how we proceed, if we are given an inequality of the form \eqref{eq:Baker-Davenport}
and an upper bound $M$ for solutions with $u\leq M$.
We start with the smallest denominator $q=q_j$ of the $j$-th convergent $\frac {p_j}{q_j}$ of $\gamma$ that exceeds $6M$. If $\varepsilon = \|\mu q\| - M \| \gamma q\|>0$,
we compute the respective upper bound $w < \frac{\log (Aq/ \varepsilon )}{\log B}$. If we get a negative $\varepsilon$, we consider the denominator 
$q_{j+1}$ of the $(j+1)$-th convergent $p_{j+1}/q_{j+1}$ instead.
If a positive $\varepsilon$ is obtained we compute the respective upper bound for $w$. If also the denominator  $q_{j+1}$ of the $(j+1)$-th convergent yields a negative $\varepsilon$
we consider the denominator of the next convergent until we obtain a positive $\varepsilon$. Let us note that it is very unlikely that after several iterations no instance occurs
with a positive $\varepsilon$, without any
good reason. Usually this reason is a rational linear dependence on $1, \gamma$ and $\mu$. If we find such a linear relation involving $1$, $\gamma$ and $\mu$, inequality \eqref{eq:Baker-Davenport}
turns into an inequality of the form
$$ 0 < \left| u' \gamma - v' \right| < AB^{-w} $$
and we are reduced to a classical approximation problem and may use the theory of continued fractions. We will treat such cases separately.  
\end{remark}


\section{Set up}

During the proof of both theorems we use the Binet formula for the Fibonacci sequence in the following form:
\begin{equation} \label{Binet}
F_k = \frac{\alpha^k-\beta^k}{\alpha - \beta} \qquad \forall k \geq 0,
\end{equation}
where $\alpha=\frac{1+\sqrt{5}}{2}$ and $\beta =\frac{1-\sqrt{5}}{2}$ are the roots of the characteristic polynomial $x^2-x-1$. Moreover, we have the inequality
\begin{align} \label{Fib_ineq}
\alpha^{k-2} \leq F_k \leq \alpha^{k-1} \qquad \forall k \geq 1.
\end{align}

Without loss of generality, we may assume that $n_1 \geq n_2 \geq 0 $ and $a_1 \geq a_2 \geq a_3 \geq 0$. Similarly, we may assume that
$m_1 \geq m_2 \geq m_3 \geq  0 $ and $t_1 \geq t_2  \geq 0$ when solving equation \eqref{eqt3A_2}.

\subsection{ Scenario for equation \eqref{eqt3A_1}}

Recall that we would like to solve 
\begin{equation*}  
F_{n_1} + F_{n_2} = 2^{a_1} + 2^{a_2} + 2^{a_3}
\end{equation*}
for $n_1, n_2, a_1, a_2$ and $a_3$. Thus we get
\begin{equation}\label{ineq1}
\alpha^{n_1 - 2} \leq F_{n_1} \leq F_{n_1} + F_{n_2} = 2^{a_1} + 2^{a_2} + 2^{a_3} \leq 3 \cdot 2^{a_1},
\end{equation}
and
\begin{equation}\label{ineq2}
2 \alpha^{n_1 -1} \geq 2F_{n_1} \geq F_{n_1} + F_{n_2} =  2^{a_1} + 2^{a_2} + 2^{a_3} > 2^{a_1}. 
\end{equation}
Hence
\begin{equation}\label{computer}
n_1 - 2 \leq a_1 \cdot \frac{\log 2}{\log \alpha} + \frac{\log 3}{\log \alpha} \quad \mbox{and} \quad n_1 - 1 \geq (a_1 -1)\cdot \frac{\log 2}{\log \alpha} ,
\end{equation}
where $\frac{\log 2}{\log \alpha} = 1.4404 \ldots$. In particular, we have $n_1 > a_1$.

In a first step, we solve equation \eqref{eqt3A_1} for all $n_1 < 360$. Inequality \eqref{computer} implies that in this case we have $a_1 < 251$. By a brute force computer
enumeration for $0 \leq n_2 \leq n_1 < 360$ and $0 \leq a_3 \leq a_2 \leq  a_1 < 251$
we found all solutions listed in the Appendix.

\subsection{ Scenario for equation \eqref{eqt3A_2}}
Recall that we would like to solve  
\begin{equation*}  
F_{m_1} + F_{m_2} + F_{m_3} =2^{t_1} + 2^{t_2} 
\end{equation*}
for $m_1, m_2, m_3, t_1$ and $t_2$. Similarly as above we obtain
\begin{equation} \label{ineq1_2}
\alpha^{m_1 - 2} \leq  F_{m_1} \leq  F_{m_1} + F_{m_2} + F_{m_3} = 2^{t_1} + 2^{t_2}\leq 2^{t_1+1}
\end{equation}
and
\begin{equation}\label{ineq2_2}
3 \alpha^{m_1 -1} \geq 3F_{m_1} \geq F_{m_1} + F_{m_2} + F_{m_3} =  2^{t_1} + 2^{t_2} > 2^{t_1}. 
\end{equation}
Thus
\begin{equation}\label{computer_2}
m_1 - 2 \leq t_1 \cdot  \frac{\log 2}{\log \alpha} + \frac{\log 2}{\log \alpha} \quad \mbox{and} \quad m_1 - 1 >  t_1 \cdot  \frac{\log 2}{\log \alpha} -\frac{\log 3}{\log \alpha}.
\end{equation}
In particular, we have $m_1 > t_1$.

We solve equation \eqref{eqt3A_2} for $0 \leq m_3 \leq m_2 \leq m_1 < 360$ and $0 \leq t_2 \leq t_1 < 251 $ by a brute force computer
enumeration and find all solutions listed in the Appendix.

By these computer searches we may assume now that $n_1 \geq 360$ for solving equation \eqref{eqt3A_1} (respectively $m_1 \geq 360$ for solving equation \eqref{eqt3A_2}). 
Moreover, we want to emphasize that the second inequality of \eqref{computer} (respectively \eqref{computer_2}) implies that $n_1 > a_1$ (respectively $m_1 > t_1$). 


\section{A first upper bound - Application of Linear forms in logarithms}

In this section, we shall establish the following two propositions concerning Diophantine equations \eqref{eqt3A_1} and \eqref{eqt3A_2} respectively.

\begin{proposition}\label{prop:bound3A_1}
 Assume that $(n_1, n_2, a_1,a_2,a_3)$ is a solution to equation \eqref{eqt3A_1} with $n_1 \geq n_2 \geq 0$ and $a_1 \geq a_2 \geq a_3 \geq 0$. Then we have that 
  $n_1 < 4.1 \cdot 10^{62}$.
\end{proposition}

\begin{proposition}\label{prop:bound3A_2}
 Assume that $(m_1, m_2, m_3,t_1,t_2)$ is a solution to equation \eqref{eqt3A_2} with $m_1 \geq m_2 \geq m_3 \geq 0$ and $t_1 \geq t_2 \geq 0$. Then we have that 
  $m_1 < 4.2 \cdot 10^{62}$.
\end{proposition}

\subsection{Proof of Proposition \ref{prop:bound3A_1}}

We follow the steps explained in the introduction. We start with\\

\noindent \textbf{Step 1:} \textit{Show that  $$\min \{ (a_1-a_2)\log 2,(n_1-n_2)\log\alpha \} < 2.61 \cdot 10^{13} \log n_1.$$}

Equation \eqref{eqt3A_1} can be rewritten as
$$ \frac{\alpha^{n_1} - \beta^{n_1}}{\sqrt{5}} + \frac{\alpha^{n_2} - \beta^{n_2}}{\sqrt{5}} =
2^{a_1} + 2^{a_2} + 2^{a_3}.$$
In the first step we consider $n_1$ and $a_1$ to be large and by collecting ``large'' terms to the left hand side of the equation, we obtain
\begin{align*}
\left| \frac{\alpha^{n_1}}{\sqrt{5}} - 2^{a_1} \right| & =  \left| 2^{a_2} + 2^{a_3}+ \frac{\beta^{n_1}}{\sqrt{5}}
- \frac{\alpha^{n_2} - \beta^{n_2}}{\sqrt{5}} \right|\\
& < 2^{a_2+1} +  \frac{\alpha^{n_2}}{\sqrt{5}} + 0.45 \\
& < 2.9 \max \{ 2^{a_2}, \alpha^{n_2} \}.
\end{align*}
Dividing through $2^{a_1}$ we get
\begin{align*}
\left| \frac{\alpha^{n_1}}{{\sqrt{5}}} 2^{-a_1}  - 1 \right| 
& < \max \left\{ 2.9 \cdot 2^{a_2 - a_1} ,\frac{2.9 \alpha^{n_2}}{2^{a_1}}  \right\}
 < \max \left\{ 2.9 \cdot 2^{a_2 - a_1} ,\frac{8.7\alpha^{n_2}}{\alpha^{n_1-2}}  \right\}. 
\end{align*}
Hence we obtain the inequality
\begin{equation}\label{Case0}
\left| \frac{\alpha^{n_1}}{{\sqrt{5}}} 2^{-a_1}  - 1 \right| 
<  22.78 \max \left\{ 2^{a_2 - a_1} , \alpha^{n_2 - n_1} \right\}.
\end{equation}

In Step 1 we consider the linear form
$${\it\Lambda} = n_1 \log \alpha - a_1 \log 2 - \log \sqrt{5}$$
and assume that $|{\it\Lambda}| \leq 0.5$.
Further, we put
$${\it\Phi} = e^{\it\Lambda}-1 = \alpha^{n_1} 2^{-a_1} {\sqrt{5}}^{-1} - 1$$
and use the theorem of Baker and W\"ustholz (Theorem \ref{BaWu}) with the data
\begin{gather*}
 \alpha_1 = \alpha, \qquad \alpha_2 = 2, \qquad \alpha_3 = \sqrt{5},\\
 b_1 = n_1, \qquad b_2 = -a_1,  \qquad b_3 = -1.
\end{gather*}
Since $n_1>a_1$ we have $B=n_1$. By simple computations, we obtain $h'(\alpha_1) = \frac{1}{2}$,
$h'(\alpha_2) = \log 2$ and $h'(\alpha_3) =  \log \sqrt{5}$.

Before we can apply Theorem \ref{BaWu} we have to show that ${\it\Phi} \neq 0$. Assume to the contrary that ${\it\Phi} = 0$,
then $\alpha^{n_1} = \sqrt{5} \cdot 2^{a_1} $. Let $\sigma \neq \mbox{id}$ be the unique non-trivial $\bbQ$-automorphism over $\bbQ(\sqrt{5})$. Then we get
$$ \alpha^{n_1} = \sqrt{5} \cdot 2^{a_1} = -\sigma\left(\sqrt{5}\cdot 2^{a_1}\right) = -\sigma\left(\alpha^{n_1}\right) = -\beta^{n_1}.$$
However, the absolute value of $\alpha^{n_1}$ is at least $\alpha^{360} > 2$ whereas the absolute value of $-\beta^{n_1}$ is at most $|\beta|^{360} < 1$. 
By this obvious contradiction we conclude that ${\it\Phi} \neq 0$.

Theorem \ref{BaWu} yields
\begin{equation*}
\log |{\it\Phi}| \geq -C(3,2) \left( \frac{1}{2} \right)
\left( \log 2 \right) \left( \log \sqrt{5} \right) \log n_1 - \log 2
\end{equation*}
and together with inequality \eqref{Case0} we have
\[\min \{ (a_1-a_2)\log 2, (n_1-n_2)\log\alpha \} < 2.61 \cdot 10^{13} \log n_1. \]
Thus we have proved so far:

\begin{lemma}\label{lem:Case0}
 Assume that $(n_1, n_2, a_1, a_2, a_3)$ is a solution to equation \eqref{eqt3A_1} with $n_1 \geq n_2 \geq 0$ and $a_1 \geq a_2 \geq a_3 \geq 0$. Then we have
 $$\min \{ (a_1-a_2)\log 2,(n_1-n_2)\log\alpha \} < 2.61 \cdot 10^{13} \log n_1.$$
\end{lemma}

Note that in the case that $|{\it\Lambda}| > 0.5$, inequality \eqref{Case0} is possible only if
either $a_1-a_2 \leq 5$ or $n_1-n_2 \leq 7$, which are covered by the bound provided by Lemma \ref{lem:Case0}.

Now we have to distinguish between
\begin{description}
\item[Case 1] $\min \{ (a_1-a_2)\log 2,(n_1-n_2)\log\alpha \}  = (a_1-a_2)\log 2$ \quad and 
\item[Case 2] $\min \{ (a_1-a_2)\log 2,(n_1-n_2)\log\alpha \}  = (n_1-n_2)\log\alpha$.
\end{description}
We will deal with these two cases in the following steps.\\

\noindent \textbf{Step 2:} \textit{We consider Case 1 and show that under the assumption that $(a_1-a_2)\log 2< 2.61 \cdot 10^{13} \log n_1$ we obtain
$$\min \left\{ (a_1 - a_3) \log 2, (n_1 - n_2) \log \alpha \right\} < 8.5 \cdot 10^{26} (\log n_1 )^2.$$}

Since we consider Case 1 we assume that
$$\min \{ (a_1-a_2)\log 2,(n_1-n_2)\log\alpha \}  = (a_1-a_2)\log 2< 2.61 \cdot 10^{13} \log n_1.$$
By collecting ``large'' terms, i.e. terms involving $n_1$, $a_1$ and $a_2$, on the left hand side, we rewrite equation \eqref{eqt3A_1} as
$$
\left| \frac{\alpha^{n_1} }{\sqrt{5}} - 2^{a_1} - 2^{a_2} \right| 
= \left| 2^{a_3} + \frac{\beta^{n_1}}{\sqrt{5}} - \frac{\alpha^{n_2}}{\sqrt{5}} + \frac{\beta^{n_2}}{\sqrt{5}} \right|
 < 2^{a_3} + \frac{\alpha^{n_2}}{\sqrt{5}} + 0.45 
$$
and obtain that
$$ \left| \frac{\alpha^{n_1} }{\sqrt{5}} - 2^{a_2} \left(2^{ a_1-a_2 } + 1 \right) \right| 
<  1.9 \max \left\{2^{a_3}, \alpha^{n_2} \right\}. $$
Dividing through $\frac{\alpha^{n_1} }{\sqrt{5}}$ we get by using inequality \eqref{ineq2}
\begin{align*}
\left| \alpha^{-n_1} 2^{a_2}  \sqrt{5} \left(2^{ a_1-a_2 } + 1 \right) - 1 \right| & < \max \left\{ \frac{1.9 \sqrt{5}}{\alpha^{n_1}} \cdot  2^{a_3}, 1.9 \sqrt{5} \alpha^{n_2-n_1} \right\}\\
& \leq \max \left\{ \frac{1.9 \sqrt{5}}{2^{a_1 - 1} \alpha } \cdot  2^{a_3}, 1.9 \sqrt{5} \alpha^{n_2-n_1} \right\}
\end{align*}
and obtain the inequality  
\begin{equation} \label{Case1}
\left| \alpha^{-n_1} 2^{a_2}   \sqrt{5} \left(2^{ a_1-a_2 } + 1 \right) - 1 \right|  < 5.26 \max \left\{ 2^{a_3 - a_1}, \alpha^{n_2-n_1} \right\}.
\end{equation}

We shall apply Theorem \ref{BaWu} to inequality \eqref{Case1}. Therefore we consider the following linear form in logarithms:
$${\it\Lambda}_1 = - n_1 \log \alpha + a_2 \log 2 + \log \left( \sqrt{5} \left(2^{ a_1-a_2 } + 1 \right) \right). $$
Let us assume for the moment that
$|{\it\Lambda}_1| \leq 0.5$. Further, we put
$${\it\Phi}_1 = e^{{\it\Lambda}_1}-1= \alpha^{-n_1} 2^{a_2} \sqrt{5} \left(2^{ a_1-a_2 } + 1 \right)  - 1 $$
and aim to apply Theorem \ref{BaWu} by taking 
\begin{gather*} \alpha_1 =  \alpha, \qquad \alpha_2 = 2,\qquad \alpha_3=\sqrt{5} \left(2^{ a_1-a_2 } + 1 \right),\\
 b_1 = -n_1,  \qquad b_2=a_2,  \qquad b_3 = 1.
\end{gather*}
Note that since $n_1>a_1>a_2$ we have $B=n_1$.
Next, we estimate the height of $\alpha_3$ by using the properties of heights and Lemma \ref{lem:Case0}:
\begin{align*}
h_0(\alpha_3) 
& \leq h_0(\sqrt{5}) + (a_1 - a_2) h_0(2) + \log 2 \\
& \leq \log \sqrt{5} + (a_1 - a_2) \log 2 + \log 2 \\
& <  2.62 \cdot 10^{13} \log n_1,
\end{align*}
which gives
$ h'(\alpha_3) <  2.62 \cdot 10^{13} \log n_1$. 
As before we have $h'(\alpha_1) = \frac{1}{2} $ and $h'(\alpha_2) = \log 2$. By a similar argument as in Step 1 we conclude that ${\it\Phi}_1 \neq 0$. 
Now, we are ready to apply Theorem \ref{BaWu} and get
\begin{align*}
\log |{\it\Phi}_1|  > & -C(3,2)  \left(  \frac{1}{2}  \right) \left(  \log 2 \right) \left( 2.62 \cdot 10^{13} \log n_1 \right)\log n_1 - \log 2\\
> & -8.49 \cdot 10^{26} (\log n_1)^2.
\end{align*}
Combining this inequality with inequality \eqref{Case1}, we obtain
\begin{equation} \label{CaseAB}
\min \left\{ (a_1 - a_3) \log 2, (n_1 - n_2) \log \alpha \right\} < 8.5 \cdot 10^{26} (\log n_1 )^2.
\end{equation}
Note that in the case that $|{\it\Lambda}_1| > 0.5$ inequality \eqref{Case1} is possible, only if either $a_1-a_3 \leq 3$ or $n_1 - n_2 \leq 4$. 
Both cases are covered by the bound provided by inequality \eqref{CaseAB}.

At this stage, we have to consider two further sub-cases.
\begin{description}
\item[Case 1A] $\min \{ (a_1-a_3)\log 2,(n_1-n_2)\log\alpha \}  = (a_1-a_3)\log 2$ \quad and 
\item[Case 1B] $\min \{ (a_1-a_3)\log 2,(n_1-n_2)\log\alpha \}  = (n_1-n_2)\log\alpha$.
\end{description}
We will deal with Case 1A in Step 3 and with Case 1B in Step 4.\\

\noindent \textbf{Step 3:} \textit{We consider Case 1A and show that under the assumption that $(a_1-a_3)\log 2< 8.5 \cdot 10^{26} (\log n_1)^2$
and $(a_1-a_2)\log 2 < 2.61 \cdot 10^{13} \log n_1$ we obtain that
$$(n_1 - n_2) \log \alpha < 2.77 \cdot 10^{40} (\log n_1 )^3.$$}

In this step we consider $n_1,a_1,a_2$ and $a_3$ to be large. By collecting ``large'' terms on the left hand side we rewrite equation \eqref{eqt3A_1} as
\begin{align*}
\left| \frac{\alpha^{n_1} }{\sqrt{5}} - 2^{a_1} - 2^{a_2} -2^{a_3}  \right| 
& = \left| \frac{\beta^{n_1}}{\sqrt{5}} - \frac{\alpha^{n_2}}{\sqrt{5}} + \frac{\beta^{n_2}}{\sqrt{5}} \right|
<  \frac{\alpha^{n_2}}{\sqrt{5}} + 0.45 
\end{align*}
and obtain that
$$ \left| \frac{\alpha^{n_1} }{\sqrt{5}} - 2^{a_1} \left( 1 + 2^{a_2-a_1} + 2^{a_3-a_1}  \right) \right| 
< 0.9 \alpha^{n_2} . $$
Dividing through $\frac{\alpha^{n_1} }{\sqrt{5}}$ yields the inequality
\begin{equation} \label{CaseA}
\left| \alpha^{-n_1} 2^{a_1} \sqrt{5} \left( 1 + 2^{a_2-a_1} + 2^{a_3-a_1}  \right) - 1 \right| 
< 2.02 \alpha^{n_2-n_1} .
\end{equation}
We want to apply Theorem \ref{BaWu} to inequality \eqref{CaseA} and consider the linear form
$${\it\Lambda}_A = - n_1 \log \alpha + a_1 \log 2 + \log \left( \sqrt{5} \left( 1 + 2^{a_2-a_1} + 2^{a_3-a_1}  \right) \right). $$
Let us assume that
$|{\it\Lambda}_A| \leq 0.5$. Further, we put
$${\it\Phi}_A = e^{{\it\Lambda}_A}-1= \alpha^{-n_1} 2^{a_1} \sqrt{5} \left( 1 + 2^{a_2-a_1} + 2^{a_3-a_1}  \right) - 1$$
and aim to apply Theorem \ref{BaWu} with
\begin{gather*}
  \alpha_1 =  \alpha, \qquad \alpha_2 = 2, \qquad \alpha_3= \sqrt{5} \left( 1 + 2^{a_2-a_1} + 2^{a_3-a_1}  \right),\\ 
  b_1 = -n_1, \qquad b_2=a_1, \qquad b_3 = 1.
\end{gather*}
 Similarly as before we get that $B = n_1$.
Next, let us estimate the height of $\alpha_3$. Using the properties of heights, Lemma \ref{lem:Case0} and inequality \eqref{CaseAB} we get
\begin{align*}
h_0(\alpha_3) 
& \leq h_0(\sqrt{5}) + (a_1 - a_2) h_0(2) + (a_1 - a_3) h_0(2) + \log 2 \\
& \leq \log \sqrt{5} + (a_1 - a_2) \log 2 + (a_1 - a_3) \log 2 + \log 2 \\
& < 2.61 \cdot 10^{13} \log n_1 + 8.5 \cdot 10^{26} (\log n_1)^2 + \log 2 \sqrt{5} \\
& <  8.51 \cdot 10^{26} (\log n_1)^2,
\end{align*}
which gives
$ h'(\alpha_3) < 8.51 \cdot 10^{26} (\log n_1)^2$. 
As before we have $h'(\alpha_1) = \frac{1}{2} $, $h'(\alpha_2) = \log 2$ and ${\it\Phi}_A \neq 0$. 
An application of Theorem \ref{BaWu} yields
\begin{align*}
\log |{\it\Phi}_A|  > & -C(3,2)  \left(  \frac{1}{2}  \right) \left(  \log 2 \right) \left( 8.51 \cdot 10^{26} (\log n_1)^2 \right)
\log n_1 - \log 2\\
> & -2.76 \cdot 10^{40} (\log n_1)^3.
\end{align*}
Combining this inequality with inequality \eqref{CaseA} we obtain
\begin{equation} \label{CaseA_bound}
(n_1 - n_2) \log \alpha < 2.77 \cdot 10^{40} (\log n_1 )^3.
\end{equation}
Note that in the case that $|{\it\Lambda}_A| > 0.5$ inequality \eqref{CaseA} is possible only if $n_1 - n_2 \leq 2$. 
This is covered by the bound provided by inequality \eqref{CaseA_bound}.\\

\noindent \textbf{Step 4:} \textit{We consider Case 1B and show that under the assumption that $(n_1-n_2)\log 2< 8.5 \cdot 10^{26} (\log n_1)^2$ 
and $(a_1-a_2)\log 2 < 2.61 \cdot 10^{13} \log n_1$
we obtain that
$$(a_1 - a_3) \log 2 < 1.39 \cdot 10^{40} (\log n_1 )^3.$$}


By collecting ``large'' terms to the left hand side, where we consider $n_1,n_2,a_1$ and $a_2$ to be large, we rewrite equation \eqref{eqt3A_1} as
\begin{align*}
\left| \frac{\alpha^{n_1} }{\sqrt{5}} + \frac{\alpha^{n_2}}{\sqrt{5}} - 2^{a_1} - 2^{a_2} \right| 
& = \left| 2^{a_3} + \frac{\beta^{n_1}}{\sqrt{5}}  + \frac{\beta^{n_2}}{\sqrt{5}} \right|
 <  2^{a_3} + 0.45 
\end{align*}
and obtain that
$$ \left| \frac{\alpha^{n_2} }{\sqrt{5}} \left( \alpha^{n_1 - n_2} + 1 \right) - 2^{a_2} \left(2^{a_1-a_2} + 1 \right) \right| 
< 1.45 \cdot 2^{a_3} . $$
Dividing through $2^{a_2} \left(2^{a_1-a_2} + 1 \right)$ we obtain the inequality
\begin{equation} \label{CaseB}
\left| \alpha^{n_2}  2^{-a_2} \left( \frac{  \alpha^{n_1 - n_2} + 1   }{\sqrt{5} \left(2^{a_1-a_2} + 1 \right)} \right)  - 1 \right|  
< 1.45 \cdot 2^{a_3 - a_1}
\end{equation} 
We want to apply Theorem \ref{BaWu} to inequality \eqref{CaseB}. Hence we consider the linear form
$${\it\Lambda}_B = n_2 \log \alpha - a_2 \log 2 + \log \left( \frac{  \alpha^{n_1 - n_2} + 1   }{\sqrt{5} \left(2^{a_1-a_2} + 1 \right)} \right)$$
and assume that
$|{\it\Lambda}_B| \leq 0.5$. Further, we put
$${\it\Phi}_B = e^{{\it\Lambda}_B}-1= \alpha^{n_2}  2^{-a_2} \left( \frac{  \alpha^{n_1 - n_2} + 1   }{\sqrt{5} \left(2^{a_1-a_2} + 1 \right)}  \right) - 1$$
and aim to apply Theorem \ref{BaWu} by taking 
\begin{gather*}
 \alpha_1 =  \alpha, \qquad \alpha_2 = 2, \qquad \alpha_3= \frac{  \alpha^{n_1 - n_2} + 1   }{\sqrt{5} \left(2^{a_1-a_2} + 1 \right)},\\
  b_1 = n_2 , \qquad b_2 = -a_2,  \qquad b_3 = 1.
\end{gather*}
  and get $B = n_1$ as in the steps before. 
Let us estimate the height of $\alpha_3$. Using the properties of heights, Lemma \ref{lem:Case0} and inequality \eqref{CaseAB} we get
\begin{align*}
h_0(\alpha_3) 
& \leq (n_1 - n_2) h_0(\alpha) + \log 2 + h_0(\sqrt{5}) + (a_1 - a_2) h_0(2) +  \log 2 \\
& = \frac{1}{2} (n_1 - n_2) \log \alpha + \log \sqrt{5} + (a_1 - a_2) \log 2 + 2\log 2 \\
& < \frac{1}{2} \left( 8.5 \cdot 10^{26} (\log n_1)^2 \right) +  2.61 \cdot 10^{13} \log n_1 + \log 4\sqrt{5} \\
& < 4.26 \cdot 10^{26} (\log n_1)^2,
\end{align*}
which gives
$ h'(\alpha_3) < 4.26 \cdot 10^{26} (\log n_1)^2$. 
A similar deduction as before yields $h'(\alpha_1) = \frac{1}{2} $, $h'(\alpha_2) = \log 2$ and ${\it\Phi}_B \neq 0$. 
Now, we apply Theorem \ref{BaWu} and get
\begin{align*}
\log |{\it\Phi}_B|  > & -C(3,2)  \left(  \frac{1}{2}  \right) \left(  \log 2 \right)\left( 4.26 \cdot 10^{26} (\log n_1)^2 \right) \log n_1 - \log 2\\
> & -1.38 \cdot 10^{40} (\log n_1)^3.
\end{align*} 
Combining this inequality with inequality \eqref{CaseB}, we obtain
\begin{equation} \label{CaseB_bound}
(a_1 - a_3) \log 2 < 1.39 \cdot 10^{40} (\log n_1 )^3.
\end{equation} 
Note that in the case of $|{\it\Lambda}_B| > 0.5$, inequality \eqref{CaseB} is possible only if $a_1 - a_3 \leq 1$ which is covered by the bound provided by inequality \eqref{CaseB_bound}.\\

\noindent \textbf{Step 5:} \textit{We consider Case 2 and show that under the assumption that $(n_1-n_2)\log \alpha< 2.61 \cdot 10^{13} \log n_1$ we obtain
$$(a_1 - a_2) \log 2 < 4.26 \cdot 10^{26} (\log n_1 )^2.$$}

Since we consider Case 2 we assume that
$$\min \{ (a_1-a_2)\log 2,(n_1-n_2)\log\alpha \}  = (n_1-n_2)\log\alpha < 2.61 \cdot 10^{13} \log n_1.$$
In this step we consider $n_1,n_2$ and $a_1$ to be large and by collecting ``large'' terms to the left hand side, we rewrite equation \eqref{eqt3A_1} as
\begin{align*}
\left| \frac{\alpha^{n_1} }{\sqrt{5}} + \frac{\alpha^{n_2}}{\sqrt{5}} - 2^{a_1} \right| 
& = \left| 2^{a_2} + 2^{a_3} + \frac{\beta^{n_1}}{\sqrt{5}}  + \frac{\beta^{n_2}}{\sqrt{5}} \right|
< 2 \cdot 2^{a_2} + 0.45  
\end{align*}
and obtain that
$$ \left| \frac{\alpha^{n_2} }{\sqrt{5}} \left( \alpha^{n_1 - n_2} + 1 \right) - 2^{a_1} \right| 
< 2.45 \cdot 2^{a_2}. $$
Dividing through $2^{a_1}$ we get the inequality 
\begin{equation}  \label{Case2}
\left|\alpha^{n_2} 2^{-a_1}  \left( \frac{  \alpha^{n_1 - n_2} + 1 }{\sqrt{5}} \right) - 1 \right| 
< 2.45 \cdot 2^{-(a_1 - a_2)}.
\end{equation}
Similarly as above we shall apply Theorem \ref{BaWu} to inequality \eqref{Case2}. Hence we consider the linear form
$${\it\Lambda}_2 =  n_2 \log \alpha - a_1 \log 2 + \log \left(\frac{  \alpha^{n_1 - n_2} + 1 }{\sqrt{5}} \right)$$
and assume that
$|{\it\Lambda}_2| \leq 0.5$. Further, we put
$${\it\Phi}_2 = e^{{\it\Lambda}_2}-1=  \alpha^{n_2} 2^{-a_1} \left(\frac{  \alpha^{n_1 - n_2} + 1 }{\sqrt{5}} \right) - 1$$
and 
\begin{gather*}
 \alpha_1 =  \alpha, \qquad \alpha_2 = 2, \qquad \alpha_3= \frac{  \alpha^{n_1 - n_2} + 1 }{\sqrt{5}} ,\\
 b_1 = n_2 , \qquad b_2 = -a_1, \qquad b_3 = 1.
\end{gather*}
  Once again this choice yields $B = n_1$. Next, let us estimate the height of $\alpha_3$. Using the properties of heights and Lemma \ref{lem:Case0} we find 
\begin{align*}
h_0(\alpha_3) 
& \leq (n_1 - n_2) h_0(\alpha) + \log 2 + h_0(\sqrt{5}) \\
& = \frac{1}{2} (n_1 - n_2) \log \alpha + \log 2 + \log \sqrt{5}  \\
& < \frac{1}{2} \left( 2.61 \cdot 10^{13} \log n_1 \right) + \log 2\sqrt{5} \\
& < 1.31 \cdot 10^{13} \log n_1,
\end{align*}
which gives
$ h'(\alpha_3) < 1.31 \cdot 10^{13} \log n_1$. 
A similar deduction as before gives  $h'(\alpha_1) = \frac{1}{2} $, $h'(\alpha_2) = \log 2$ and ${\it\Phi}_2 \neq 0$. 
Thus by applying Theorem \ref{BaWu} we get
\begin{align*}
\log |{\it\Phi}_2|  > & -C(3,2)  \left(  \frac{1}{2}  \right) \left(  \log 2 \right)\left( 1.31 \cdot 10^{13} \log n_1 \right) \log n_1 - \log 2\\
> & -4.25 \cdot 10^{26} (\log n_1)^2.
\end{align*}
Combining this inequality together with inequality \eqref{Case2}, we obtain
\begin{equation} \label{Case2_bound}
(a_1 - a_2) \log 2 < 4.26 \cdot 10^{26} (\log n_1 )^2.
\end{equation}
Note that in the case of $|{\it\Lambda}_2| > 0.5$, inequality \eqref{Case2} is possible only if $a_1 - a_2 \leq 2$ which is covered by the bound provided by inequality \eqref{Case2_bound}.\\

\noindent \textbf{Step 6:} \textit{We continue to consider Case 2 and show that under the assumption that $(n_1-n_2)\log \alpha< 2.61 \cdot 10^{13} \log n_1$ and 
$(a_1 - a_2) \log 2 < 4.26 \cdot 10^{26} (\log n_1 )^2$ we obtain
$$(a_1 - a_3) \log 2 < 1.39 \cdot 10^{40}(\log n_1 )^3.$$}

We shall apply once more Theorem \ref{BaWu} to obtain an upper bound for $(a_1 - a_3) \log 2$. The derivation is very similar to Case 1B. 
By collecting ``large'' terms on the left hand side, we rewrite equation~\eqref{eqt3A_1} as
\begin{align*}
\left| \frac{\alpha^{n_1} }{\sqrt{5}} + \frac{\alpha^{n_2}}{\sqrt{5}} - 2^{a_1} - 2^{a_2} \right| 
& = \left| 2^{a_3} + \frac{\beta^{n_1}}{\sqrt{5}}  + \frac{\beta^{n_2}}{\sqrt{5}} \right|
 <  2^{a_3} + 0.45 .
\end{align*}
By the same derivation as in Step 4 we obtain inequality \eqref{CaseB}, i.e.
\begin{equation*} 
\left| \alpha^{n_2}  2^{-a_2} \left( \frac{  \alpha^{n_1 - n_2} + 1   }{\sqrt{5} \left(2^{a_1-a_2} + 1 \right)}  \right)  - 1 \right|  
< 1.45 \cdot 2^{a_3 - a_1}.
\end{equation*}
We have the same setting as in Case 1B, except that the estimate for the height of $\alpha_3$ becomes
\begin{align*}
h_0(\alpha_3) 
& \leq (n_1 - n_2) h_0(\alpha) + \log 2 + h_0(\sqrt{5}) + (a_1 - a_2) h_0(2) +  \log 2 \\
& = \frac{1}{2} (n_1 - n_2) \log \alpha + \log \sqrt{5} + (a_1 - a_2) \log 2 + 2\log 2 \\
& < \frac{1}{2} \left( 2.61 \cdot 10^{13} \log n_1 \right) +  4.26 \cdot 10^{26} (\log n_1)^2 + \log 4\sqrt{5} \\
& < 4.27 \cdot 10^{26} (\log n_1)^2,
\end{align*}
which gives $ h'(\alpha_3) < 4.27 \cdot 10^{26} (\log n_1)^2$ instead of $ h'(\alpha_3) < 4.26 \cdot 10^{26} (\log n_1)^2$ . 
Therefore by applying Theorem \ref{BaWu} similarly as before we obtain 
\begin{equation} \label{Case2_bound2}
(a_1 - a_3) \log 2 < 1.39 \cdot 10^{40} (\log n_1 )^3
\end{equation}
which coincides with inequality \eqref{CaseB_bound}.
Table \ref{TableCase1log} summarizes our results obtained so far.\\ 

\begin{table}[ht] 
\caption{Summary of results} \label{TableCase1log} 
\centering 
\begin{tabular}{c c c c} 
\hline 
\textbf{Upper bound of}  & \textbf{Case 1A} &\textbf{ Case 1B} & \textbf{Case 2 }\\ [0.5ex]
\hline 
$(a_1 - a_2) \log 2$ 
& $2.61 \cdot 10^{13} \log n_1$ & $2.61 \cdot 10^{13} \log n_1$ & $4.26 \cdot 10^{26} (\log n_1)^2$ \\ 
$(a_1 - a_3) \log 2$ 
& $8.51 \cdot 10^{26} (\log n_1)^2$  & $1.39 \cdot 10^{40} (\log n_1)^3$  & $1.39 \cdot 10^{40} (\log n_1)^3$ \\
$(n_1 - n_2) \log \alpha$ 
& $2.77 \cdot 10^{40} (\log n_1)^3$  & $8.5 \cdot 10^{26} (\log n_1)^2$  & $2.61 \cdot 10^{13} \log n_1$ \\
[1ex] 
\hline 
\end{tabular}
\end{table}

\noindent \textbf{Step 7:} \textit{We assume the bounds given in Table \ref{TableCase1log} and show that
$n_1 \log \alpha < 4.54 \cdot 10^{53} (\log n_1)^4$, hence $n_1< 4.1 \cdot 10^{62}$.}\\

We have to apply Theorem \ref{BaWu} once more. This time we rewrite equation \eqref{eqt3A_1} as
\begin{align*}
\left| \frac{\alpha^{n_1} }{\sqrt{5}} \left(1 + \alpha^{n_2-n_1} \right) - 2^{a_1} (1+ 2^{a_2-a_1} + 2^{a_3-a_1} ) \right| 
& = \left| \frac{\beta^{n_1}}{\sqrt{5}}  + \frac{\beta^{n_2}}{\sqrt{5}} \right|
<  0.45 .
\end{align*}
Dividing through $\frac{\alpha^{n_1} }{\sqrt{5}} \left(1 + \alpha^{n_2-n_1} \right) $ we obtain the inequality
\begin{equation} \label{Case3}
\left| \alpha^{-n_1} 2^{a_1} \left( \frac{\sqrt{5} (1+ 2^{a_2-a_1} + 2^{a_3-a_1} )}  { 1 + \alpha^{n_2-n_1}  }  \right) -1 \right| 
<  1.01\alpha^{-n_1}.
\end{equation}
In this final step we consider the linear form
$${\it\Lambda}_3 = - n_1 \log \alpha + a_1 \log 2 + \log \left( \frac{\sqrt{5} (1+ 2^{a_2-a_1} + 2^{a_3-a_1} )}  { 1 + \alpha^{n_2-n_1} } \right) $$
and assume that $|{\it\Lambda}_3| \leq 0.5$. Further, we put
$${\it\Phi}_3 = e^{{\it\Lambda}_3} - 1 = \alpha^{-n_1} 2^{a_1} \left( \frac{\sqrt{5} (1+ 2^{a_2-a_1} + 2^{a_3-a_1} )}  { 1 + \alpha^{n_2-n_1}  } \right) -1.$$
We take 
\begin{gather*}
 \alpha_1 = \alpha, \qquad \alpha_2 = 2, \qquad \alpha_3 = \frac{\sqrt{5}(1+ 2^{a_2-a_1} + 2^{a_3-a_1})}{1+\alpha^{n_2-n_1}},\\
 b_1 = -n_1,\qquad  b_2 = a_1, \qquad b_3 = 1.
\end{gather*}
Thus we have $B = n_1$. By the results in Table \ref{TableCase1log} and similar computations done before we obtain
\begin{align*}
h_0 \left( \alpha_3 \right)& \leq  h_0 \left( \sqrt{5} \right) + (a_1 - a_2) \ h_0 ( 2 ) + (a_1 - a_3) \ h_0 ( 2 )  + (n_1 - n_2) h_0 ( \alpha ) + 2\log 2  \\
& \leq  (a_1 - a_2) \log 2 + (a_1 - a_3) \log 2 + \frac{1}{2}(n_1 - n_2) \log \alpha + \log 4 \sqrt{5}  \\
& <  1.4 \cdot 10^{40} (\log n_1)^3,
\end{align*}
which gives
$ h'(\alpha_3) < 1.4 \cdot 10^{40} (\log n_1)^3$. 
As before we have $h'(\alpha_1) = \frac{1}{2} $, $h'(\alpha_2) =  \log 2$ and ${\it\Phi}_3 \neq 0$. 
Now an application of Theorem \ref{BaWu} yields
$$
\log |{\it\Phi}_3|
 > -C(3,2)\left( \frac{1}{2} \right) \left( \log 2 \right) \left(1.4 \cdot 10^{40} (\log n_1)^3 \right)
   \log n_1 - \log 2.
$$
Combining this inequality with inequality~\eqref{Case3} we get
$$ n_1 \log \alpha < 4.54 \cdot 10^{53} (\log n_1)^4$$
which yields 
$$n_1 < 4.1 \cdot 10^{62}.$$
Similarly as in the cases above the assumption $|{\it\Lambda}_3| > 0.5$ leads in view of inequality \eqref{Case3} to $n_1 \leq 0$ which is impossible. Thus Proposition \ref{prop:bound3A_1} is established.

\subsection{Proof of Proposition \ref{prop:bound3A_2}}

Since the deduction of an upper bound for solutions to \eqref{eqt3A_2} is similar to the proof of Proposition \ref{prop:bound3A_1} we only sketch the argument.
In the case of equation \eqref{eqt3A_2}, we have
$$ \frac{\alpha^{m_1} - \beta^{m_1}}{\sqrt{5}} + \frac{\alpha^{m_2} - \beta^{m_2}}{\sqrt{5}} + \frac{\alpha^{m_3} - \beta^{m_3}}{\sqrt{5}} =
2^{t_1} + 2^{t_2} .$$

\noindent \textbf{Step 1:} \textit{Show that 
$$\min \{ (t_1-t_2)\log 2,(m_1-m_2)\log\alpha \} < 2.61 \cdot 10^{13} \log m_1.$$}

First, we rearrange equation \eqref{eqt3A_2} and make use of inequalities \eqref{ineq1_2} and \eqref{ineq2_2} to get
\begin{equation}\label{Case0_2}
\left| \frac{\alpha^{m_1}2^{-t_1}}{{\sqrt{5}}} - 1 \right| 
<  14.67 \max \left\{ 2^{t_2 - t_1} , \alpha^{m_2 - m_1} \right\}.
\end{equation}
We consider ${\it\Gamma} = m_1 \log \alpha - t_1 \log 2 - \log \sqrt{5}$ 
with $|{\it\Gamma}| \leq 0.5$.
Further, we put
$${\it\Psi} = e^{\it\Gamma}-1 = \alpha^{m_1} 2^{-t_1} {\sqrt{5}}^{-1} - 1$$
and apply the theorem of Baker and W\"ustholz (Theorem \ref{BaWu}) with the data
\begin{gather*}
 \alpha_1 = \alpha, \qquad  \alpha_2 = 2, \qquad  \alpha_3 = \sqrt{5},\\
 b_1 = m_1, \qquad  b_2 = -t_1,   \qquad  b_3 = -1,
\end{gather*}
 i.e.  $B=m_1$. By a simple computation, we obtain $h'(\alpha_1) = \frac{1}{2}$,
$h'(\alpha_2) = \log 2$ and $h'(\alpha_3) =  \log \sqrt{5}$. Similarly as in the proof of Proposition \ref{prop:bound3A_1} we may assume that ${\it\Psi} \neq 0$.
Then Theorem \ref{BaWu} yields
\begin{equation*}
\log |{\it\Psi}| \geq -C(3,2) \left( \frac{1}{2} \right)
\left( \log 2 \right) \left( \log \sqrt{5} \right) \log m_1 - \log 2
\end{equation*}
and together with inequality \eqref{Case0_2} we have
\[\min \{ (t_1-t_2)\log 2, (m_1-m_2)\log\alpha \} < 2.61 \cdot 10^{13} \log m_1. \]
Thus instead of Lemma \ref{lem:Case0} we obtain now

\begin{lemma}\label{lem:Case0_2}
 Assume that $(m_1, m_2, m_3, t_1, t_2)$ is a solution to equation \eqref{eqt3A_2} with $m_1 \geq m_2 \geq m_3 \geq 0$ and $t_1 \geq t_2 \geq 0$. Then we have
 $$\min \{ (t_1-t_2)\log 2,(m_1-m_2)\log\alpha \} < 2.61 \cdot 10^{13} \log m_1.$$
\end{lemma}

The scenarios for which $|{\it\Gamma}| > 0.5$ can be easily dealt with. Now we have to distinguish between two cases:
\begin{description}
\item[Case 1] $\min \{ (t_1-t_2)\log 2,(m_1-m_2)\log\alpha \}  = (m_1-m_2)\log \alpha $ \quad and 
\item[Case 2] $\min \left\{ (t_1 - t_2) \log 2, (m_1 - m_3) \log \alpha \right\}  = (t_1-t_2)\log 2$
\end{description}
We will deal with these cases in the following steps.\\

\noindent \textbf{Step 2:} \textit{We consider Case 1 and show that under the assumption that $(m_1-m_2)\log\alpha < 2.61 \cdot 10^{13} \log m_1$ we obtain
$$\min \left\{ (t_1 - t_2) \log 2, (m_1 - m_3) \log \alpha \right\} < 4.26 \cdot 10^{26} (\log m_1 )^2.$$}

We rearrange equation \eqref{eqt3A_2} and make use of inequalities \eqref{ineq1_2} and \eqref{ineq2_2} to get
\begin{equation} \label{Case1_2}
| {\it\Psi}_1 | = \left| \frac{\alpha^{m_2} 2^{-t_1} \left(\alpha^{ m_1-m_2 } + 1 \right) }{{\sqrt{5}}}   - 1  \right|  < 12.31 \max \left\{ 2^{-(t_1 - t_2)}, \alpha^{-(m_1-m_3)} \right\}.
\end{equation}
We apply Theorem \ref{BaWu} to inequality \eqref{Case1_2} by taking
$b_1=m_2$, $b_2=-t_1$ and $b_3=1$, i.e. $B = m_1$ since $m_1>m_2,t_1$. Further, we choose $\alpha_1 =  \alpha$, $\alpha_2 = 2$ and $\alpha_3= \frac{ \alpha^{ m_1-m_2 } + 1 }{{\sqrt{5}}} $.
Note that by our standard arguments we obtain that
$ h'(\alpha_3) <  1.31 \cdot 10^{13} \log m_1 $ and ${\it\Psi}_1  \neq 0$. 
Finally we get
\begin{equation*} \label{CaseAB_2}
\min \left\{ (t_1 - t_2) \log 2, (m_1 - m_3) \log \alpha \right\} < 4.26 \cdot 10^{26} (\log m_1 )^2.
\end{equation*}

At this stage, we have to consider the following two sub-cases for Case 1:
\begin{description}
\item[Case 1A] $\min \left\{ (t_1 - t_2) \log 2, (m_1 - m_3) \log \alpha \right\} =(m_1-m_3)\log \alpha$ and 
\item[Case 1B] $\min \left\{ (t_1 - t_2) \log 2, (m_1 - m_3) \log \alpha \right\} =(t_1 - t_2) \log 2$.
\end{description}
We will deal with these sub-cases in the steps below.\\

\noindent \textbf{Step 3:} \textit{We consider Case 1A and show that under the assumption that $(m_1-m_2)\log\alpha < 2.61 \cdot 10^{13} \log m_1$ 
and $(m_1 - m_3) \log \alpha < 4.26 \cdot 10^{26} (\log m_1 )^2$ we obtain that
$$(t_1 - t_2) \log 2 < 6.94 \cdot 10^{39} (\log m_1 )^3.$$}

We rearrange equation \eqref{eqt3A_2} and make use of inequalities \eqref{ineq1_2} and \eqref{ineq2_2} to get
\begin{equation} \label{CaseA_2}
| {\it\Psi}_A | = \left| \frac{\alpha^{m_1} 2^{-t_1}\left( 1 + \alpha^{m_2-m_1} + \alpha^{m_3-m_1}  \right) }{{\sqrt{5}}}   - 1 \right| 
< 1.9 \cdot 2^{t_2-t_1}.
\end{equation}
We apply Theorem \ref{BaWu} to inequality \eqref{Case1_2} with 
$B = m_1$, $ \alpha_1 =  \alpha$,  $\alpha_2 = 2$, $\alpha_3=  \frac{\left( 1 + \alpha^{m_2-m_1} + \alpha^{m_3-m_1}  \right)}{\sqrt{5}} $. Note that we have
$ h'(\alpha_3) < 2.14 \cdot 10^{26} (\log m_1)^2$ and ${\it\Psi}_A \neq 0$. Therefore, we get
\begin{equation*} 
(t_1 - t_2) \log 2 < 6.94 \cdot 10^{39} (\log m_1 )^3.
\end{equation*}

\noindent \textbf{Step 4:} \textit{We consider Case 1B and show that under the assumption that $(m_1-m_2)\log\alpha < 2.61 \cdot 10^{13} \log m_1$ 
and $(t_1 - t_2) \log 2 < 4.26 \cdot 10^{26} (\log m_1 )^2$ we obtain that
$$(m_1 - m_3) \log \alpha < 1.4 \cdot 10^{40} (\log m_1 )^3.$$}

We rearrange equation \eqref{eqt3A_2} and make use of inequalities \eqref{ineq1_2} and \eqref{ineq2_2} to get
\begin{equation} \label{CaseB_2}
| {\it\Psi}_B | = \left| \alpha^{-m_2}  2^{t_2} \sqrt{5}\left( \frac{ 2^{t_1 - t_2} + 1   }{ \alpha^{m_1-m_2} + 1 }  \right)  - 1 \right|  
< 3.02 \alpha^{m_3 - m_1}
\end{equation}
We apply Theorem \ref{BaWu} to inequality \eqref{CaseB_2} by taking
$ B = m_1$, $\alpha_1 =  \alpha$, $\alpha_2 = 2$ and $\alpha_3= \frac{\sqrt{5}  \left(  2^{t_1 - t_2} + 1   \right) }{ \alpha^{m_1-m_2} + 1 }$. With this choice we have 
$ h'(\alpha_3) < 4.27 \cdot 10^{26} (\log n_1)^2$ and  ${\it\Psi}_B \neq 0$.
and we obtain
\begin{equation*} \label{CaseB_2_bound}
(m_1 - m_3) \log \alpha < 1.4 \cdot 10^{40} (\log m_1 )^3.
\end{equation*}

\noindent \textbf{Step 5:} \textit{We consider Case 2 and show that under the assumption that $(t_1-t_2)\log 2 < 2.61 \cdot 10^{13} \log m_1$ we obtain
$$(m_1 - m_2) \log \alpha < 8.5 \cdot 10^{26} (\log m_1 )^2.$$}

We rearrange equation \eqref{eqt3A_2} and make use of inequalities \eqref{ineq1_2} and \eqref{ineq2_2} to get
\begin{equation}  \label{Case2_2}
| {\it\Psi}_2 | = \left|\alpha^{-m_2} 2^{t_2} \sqrt{5}\left( 2^{t_1-t_2} +1 \right) - 1 \right| 
< 4.03 \alpha^{-(m_1 - m_2)}.
\end{equation}
We apply Theorem \ref{BaWu} to inequality \eqref{Case2_2} by taking
$B = m_1$, $\alpha_1 =  \alpha$, $\alpha_2 = 2$, $\alpha_3= \sqrt{5}\left( 2^{t_1-t_2} +1 \right)$. In this case we have that
$ h'(\alpha_3) < 2.62 \cdot 10^{13} \log m_1 $ and also ${\it\Psi}_2 \neq 0$. 
Therefore we get
\begin{equation*} \label{Case2_2_bound}
(m_1 - m_2) \log \alpha < 8.5 \cdot 10^{26} (\log m_1 )^2.
\end{equation*}

\noindent \textbf{Step 6:} \textit{We continue to consider Case 2 and show that under the assumption that $(t_1-t_2)\log 2 < 2.61 \cdot 10^{13} \log m_1$ 
and $(m_1 - m_2) \log \alpha < 8.5 \cdot 10^{26} (\log m_1 )^2$ we obtain that
$$(m_1 - m_3) \log \alpha < 1.38 \cdot 10^{40} (\log m_1 )^3.$$}

Again we apply Theorem \ref{BaWu} to obtain an upper bound for $(m_1 - m_3) \log \alpha$. The derivation is very similar to Case 1B. In particular, we have
$$
\left| \alpha^{-m_2}  2^{t_2} \sqrt{5}\left( \frac{ 2^{t_1 - t_2} + 1   }{ \alpha^{m_1-m_2} + 1 }  \right)  - 1 \right|  
< 3.02 \alpha^{m_3 - m_1}
$$
and the same setting as in Case 1B, except that 
$h'(\alpha_3) < 4.26 \cdot 10^{26} (\log m_1)^2$. 
Therefore Theorem~\ref{BaWu} gives us
\begin{equation*} 
(m_1 - m_3) \log \alpha < 1.38 \cdot 10^{40} (\log m_1 )^3.
\end{equation*}

Table \ref{TableCase2log} summarizes our results obtained so far.\\

\begin{table}[ht] 
\caption{Summary of results} \label{TableCase2log}
\centering 
\begin{tabular}{c c c c} 
\hline 
\textbf{Upper bound of}  & \textbf{Case 1A} &\textbf{ Case 1B} & \textbf{Case 2 }\\ [0.5ex] 
\hline  line
$(m_1 - m_2) \log \alpha $ 
& $2.61 \cdot 10^{13} \log m_1$ & $2.61 \cdot 10^{13} \log m_1$ & $8.5 \cdot 10^{26} (\log m_1)^2$ \\ 
$(m_1 - m_3) \log \alpha $ 
& $4.26 \cdot 10^{26} (\log m_1)^2$  & $1.4 \cdot 10^{40} (\log m_1)^3$  & $1.38 \cdot 10^{40} (\log m_1)^3$ \\
$(t_1 - t_2) \log 2$ 
& $6.94 \cdot 10^{39} (\log m_1)^3$  & $4.26 \cdot 10^{26} (\log m_1)^2$  & $2.61 \cdot 10^{13} \log m_1$ \\
[1ex] 
\hline 
\end{tabular}
\end{table}

\noindent \textbf{Step 7:} \textit{We assume the bounds given in Table \ref{TableCase2log} and show that $m_1< 4.2 \cdot 10^{62}$.}\\

Once again we have to apply Theorem \ref{BaWu}. We rearrange equation \eqref{eqt3A_2} and make use of inequalities~\eqref{ineq1_2} and \eqref{ineq2_2} to get
\begin{equation} \label{Case3_2}
| {\it\Psi}_3 | = \left| \alpha^{-m_1} 2^{t_1} \left( \frac{\sqrt{5} (1+ 2^{t_2-t_1}  )}  { 1 + \alpha^{m_2-m_1} + \alpha^{m_3-m_1} }  \right) -1 \right| 
<  2.02 \alpha^{-m_1}.
\end{equation}
In our last step we apply Theorem \ref{BaWu} to inequality \eqref{Case3_2} by taking
$ B=m_1$, $\alpha_1 = \alpha$, $\alpha_2 = 2$, $\alpha_3 = \frac{\sqrt{5} (1+ 2^{t_2-t_1}  )}  { 1 + \alpha^{m_2-m_1} + \alpha^{m_3-m_1} }$. By our usual arguments we show that
$h'(\alpha_3) < 1.41 \cdot 10^{40} (\log m_1)^3$ and ${\it\Psi}_3 \neq 0$.
Thus we get
$$m_1 < 4.2 \cdot 10^{62},$$
hence Proposition \ref{prop:bound3A_2} is established.

\begin{remark}
The theorem of Baker and W\"ustholz (cf. Theorem \ref{BaWu}) \cite{bawu93} has a significant role in the development of linear forms in logarithms. 
The final structure for the lower bound for linear forms in logarithms without an explicit determination of the constant involved has been established by W\"ustholz \cite{wu88} and the precise determination of that constant is the central aspect of \cite{bawu93}
(see also \cite{bawu07}). 
The reader may note that slightly sharper bounds for $n_1$ and $m_1$ could be obtained by using Matveev's result \cite{matveev00} instead. 
However, the improvement is insignificant in view of our next step, i.e. the use of the method of Baker and Davenport (Lemma \ref{Reduction}), in which our upper bounds for $n_1$ and $m_1$
are further reduced to a great extent. 
\end{remark}


\section{Reduction of the bound}
In our final step we reduce the huge upper bound for $n_1$ obtained in Proposition \ref{prop:bound3A_1} (respectively $m_1$ in Proposition \ref{prop:bound3A_2})
by applying several times Lemma \ref{Reduction}.

\subsection{Proof of Theorem \ref{Main3A_1}}

First, we consider inequality \eqref{Case0} and recall that
$${\it\Lambda} = n_1 \log \alpha - a_1 \log 2 - \log \sqrt{5}.$$
For technical reasons we assume that $\min\{n_1-n_2, a_1-a_2, a_1 - a_3\} \geq 20$. In the case that
this condition fails we do the following:
\begin{itemize}
\item if $a_1-a_2 < 20$ but $a_1-a_3, n_1-n_2 \geq 20$, we consider inequality \eqref{Case1}, i.e. we go to Step 2;
\item if $a_1-a_2, a_1-a_3 < 20$ but $n_1-n_2 \geq 20$, we consider inequality \eqref{CaseA}, i.e. we go to Step 3;
\item if $a_1-a_2, n_1-n_2 < 20$ but $a_1-a_3 \geq 20$, we consider inequality \eqref{CaseB}, i.e. we go to Step 4;
\item if $n_1-n_2 < 20$ but $a_1-a_2, a_1-a_3 \geq 20$, we consider inequality \eqref{Case2}, i.e. we go to Step~5; then we consider inequality \eqref{CaseB}, i.e. we go to Step 6;
\item if all $a_1-a_2, a_1-a_3, n_1-n_2 < 20$, we consider inequality \eqref{Case3}, i.e. we go to Step 7.
\end{itemize}

\noindent \textbf{Step 1:} \textit{We show that $a_1-a_2 \leq 218$ or $n_1-n_2 \leq 315$.}\\

Let us start by considering inequality \eqref{Case0}. Since we assume that $\min\{n_1-n_2, a_1-a_2\} \geq 20$
we get $|{\it\Phi}| = |e^{\it\Lambda}-1|<\frac{1}{4}$, hence $|{\it\Lambda}|< \frac{1}{2}$.
And, since $|x| < 2|e^x-1|$ holds for all $x \in (-\frac{1}{2}, \frac{1}{2})$ we get
$ |{\it\Lambda}|<45.56 \max \{ 2^{a_2 - a_1}, \alpha^{n_2 - n_1} \} $. 
Then we have the inequality
\begin{equation*}
\begin{split}
0< \left| n_1 \cdot \frac{\log \alpha}{\log 2 } -a_1 + \frac{\log(1/\sqrt{5})}{\log 2 } \right| < &
\max\left\{  \frac{45.56}{\log 2} \cdot 2^{-(a_1-a_2)}, \frac{45.56}{\log 2} \alpha^{-(n_1-n_2)} \right\}\\
< & \max\left\{ 66 \cdot 2^{-(a_1-a_2)} , 66 \alpha^{-(n_1-n_2)} \right\}
\end{split}
\end{equation*} 
and we apply the algorithm described in Remark \ref{RemarkReduction} with
\[ \gamma = \frac{\log \alpha}{\log 2}, \quad \mu = \frac{\log(1/\sqrt{5})}{\log 2}, \quad (A,B)=(66, 2) \mbox{ or } (66, \alpha). \]
Let us be a bit more precise.
We note that $\gamma$ is irrational since $2$ and $\alpha$ are multiplicatively independent, hence Lemma \ref{Reduction} is applicable.
Let $\gamma = [s_0, s_1, s_2, \dotso]=[0,1,2,3,1,2,3,2,4,2,1,2,11,\dotso]$ be the continued fraction expansion of $\gamma$.
Moreover, we choose $M = 4.1 \cdot 10^{62}$ and consider the $125$-th convergent
\[\frac{p_{125}}{q_{125}}=
\frac{2028312018571414606476009600985599840687019168230545776285240837}{2921621381175511963618293669947470310883223581600886270426241482}, \]
with $q=q_{125} > 6M$. 
This yields $\varepsilon > 0.24$ and therefore either
\[ a_1-a_2 \leq \frac{\log(66q/0.24)}{\log 2} < 219 \quad \mbox{or} \quad n_1-n_2 \leq \frac{\log(66q/0.24)}{\log \alpha} < 316. \]
Thus, we have either $a_1-a_2 \leq 218$ or $n_1-n_2 \leq 315$.

From this result we distinguish between
\begin{description}
 \item[Case 1] $a_1-a_2 \leq 218$ \quad and
 \item[Case 2] $n_1-n_2 \leq 315$.
\end{description}

\noindent \textbf{Step 2:} \textit{We consider Case 1 and show that under the assumption that $a_1-a_2 \leq 218$ we have that 
$a_1-a_3 \leq  225$ or $n_1-n_2 \leq  324$.}\\

In this step we consider inequality~\eqref{Case1} and assume that $a_1-a_3, n_1 - n_2 \geq 20$. Recall that
$${\it\Lambda}_1 = - n_1 \log \alpha + a_2 \log 2 + \log \left( \sqrt{5} \left(2^{ a_1-a_2 } + 1 \right)\right) $$
and inequality \eqref{Case1} yields that
$ |{\it\Lambda}_1| < 10.52 \max \left\{ 2^{-(a_1 - a_3)}, \alpha^{-(n_1-n_2)} \right\}$. 
Then we get
\[ 0< \left| n_1 \cdot \frac{\log \alpha}{\log 2}  -a_2 +
\frac{\log \left(1/\left({\sqrt{5} \left(2^{ a_1-a_2 } + 1 \right)} \right)\right) }{\log 2} \right|
< 16 \max \left\{ 2^{-(a_1 - a_3)}, \alpha^{-(n_1-n_2)} \right\}.\]
We apply the algorithm explained in Remark \ref{RemarkReduction} 
again with the same $\gamma$ and $M$ as in Step 1,
but now we choose $(A,B)=(16, 2)$ or $(16, \alpha)$ and 
\[ \mu=\mu_k = \frac{\log \left(1/\left({\sqrt{5} \left(2^k + 1 \right)} \right)\right) }{\log 2}\]
for each possible value of $a_1-a_2=k=0,1,\dots, 218$. With these parameters we run our algorithm
and obtain for each instance a new and rather small upper bound either for $a_1-a_3$ or $n_1-n_2$. In particular
$$q_{128}=49310467685085622966403899548743583219671853934492723134649593651$$
is the largest denominator that appeared in applying our algorithm. Overall, we obtain
\[ a_1-a_3 \leq 225 \quad \mbox{or} \quad n_1-n_2 \leq 324. \]

Within Case 1 we have to distinguish between two further sub-cases:
\begin{description}
 \item[Case 1A] $a_1-a_3 \leq 225$ \quad and
 \item[Case 1B] $n_1-n_2 \leq 324$.
\end{description}

\noindent \textbf{Step 3:} \textit{We consider Case 1A and show that under the assumption that $a_1-a_2 \leq 218$ and $a_1-a_3\leq 225$ we have that 
$n_1-n_2 \leq  334$.}\\

In this step we consider inequality~\eqref{CaseA} and assume that $n_1 - n_2 \geq 20$. Recall that
$${\it\Lambda}_A = - n_1 \log \alpha + a_1 \log 2 + \log \left(\sqrt{5} \left( 1 + 2^{a_2-a_1} + 2^{a_3-a_1}  \right) \right) $$
and inequality \eqref{CaseA} yields that
$ |{\it\Lambda}_A| < 4.04 \alpha^{-(n_1-n_2)} $. 
Then we get
\[ 0< \left| n_1 \cdot \frac{\log \alpha}{\log 2} - a_1  + \frac{\log \left(1 / \sqrt{5} \left( 1 + 2^{a_2-a_1} + 2^{a_3-a_1}  \right) \right) }{\log 2} \right|
<  6 \alpha^{-(n_1-n_2)}.\]

We proceed as in Remark \ref{RemarkReduction} with the same $\gamma$ and $M$ as in Step 1, but we use $(A,B)=(6, \alpha)$ instead. Moreover we consider 
\[ \mu=\mu_{k,l} = \frac{\log \left(1 / \sqrt{5} \left( 1 + 2^{-k} + 2^{-l}  \right) \right) }{\log 2} \]
 for each possible value of $a_1-a_2=k=0,1,\dots, 218$ and $a_1-a_3=l=0,1,\dots, 225$ (with respect to the obvious condition that $a_1 - a_2 \leq a_1 - a_3$). 
 As in the previous step we apply the algorithm described in Remark \ref{RemarkReduction}
 to each instance $(k,l)$ and start
with the $125$-th convergent $\frac{p}{q} = \frac{p_{125}}{q_{125}}$ of $\gamma$
as before and continue with the algorithm until a positive $\varepsilon$ is obtained for every $k$ and $l$. Thus we can compute a new upper bound for $n_1 - n_2$ by the formula
$ n_1-n_2 < \frac{\log(6 q /\varepsilon)}{\log \alpha}  $
for the respective choices of $q$ and $\varepsilon$. Overall we obtain that 
\[ n_1-n_2 \leq 334. \]

\noindent \textbf{Step 4:} \textit{We consider Case 1B and show that under the assumption that $a_1-a_2 \leq 218$ and $n_1-n_2\leq 324$ we have that 
$a_1-a_3 \leq  233$.}\\

Thus we consider inequality~\eqref{CaseB} and assume that $a_1 - a_3 \geq 20$. In view of Step 6 we perform the following reduction by considering $a_1 - a_2 \leq 224$ instead of $a_1 - a_2 \leq 218$. 
Note that the same inequality \eqref{CaseB} will be used once more with a slightly higher upper bound $a_1 - a_2 \leq 224$ in Step~6.  Recall that
$${\it\Lambda}_B = n_2 \log \alpha - a_2 \log 2 + \log \left(\frac{  \alpha^{n_1 - n_2} + 1   }{\sqrt{5} \left(2^{a_1-a_2} + 1 \right)} \right) $$
and inequality \eqref{CaseB} yields that
$ |{\it\Lambda}_B| < 2.9 \cdot 2^{-(a_1-a_3)} $. 
Then we get
\[ 0< \left| n_2 \cdot \frac{\log \alpha}{\log 2} - a_2  + \frac{\log \left( (\alpha^{n_1 - n_2} + 1) / \left( \sqrt{5} \left(2^{a_1-a_2} + 1 \right) \right)\right) }{\log 2} \right|
< 5 \cdot 2^{-(a_1-a_3)} .\]
We apply our algorithm with the same $\gamma$ and $M$ as in the previous steps, but we use $(A,B)=(5, 2)$ and 
\[ \mu=\mu_{k,r} = \frac{\log \left( (\alpha^{r} + 1) / \left( \sqrt{5} \left(2^{k} + 1 \right) \right)\right) }{\log 2}\]
for each possible value of $a_1-a_2=k=0,1,\dots, 224$ and $n_1-n_2=r=0,1,\dots, 324$. 
We run our algorithm starting with $q=q_{125}$ and compute the upper bound for $a_1 - a_3$  by the formula  
$a_1-a_3 < \frac{\log(5q /\varepsilon)}{\log 2}$ 
for respective choices of $q$ and $\varepsilon$, provided the algorithm terminates. For those pairs $(k,r)$ for which the algorithm terminates we obtain  
\[  a_1-a_3  \leq 233. \] 

However, in case that $(k, r) \in \lbrace  (0, 2), (0, 6), (2, 10), (4, 18)\rbrace$ problems arise and our algorithm does not terminate. This is because in these cases there exist multiplicative
dependences between $\mu_{k,r}$, $2$ and $\alpha$. In particular, one can easily check that
$$ \frac{\alpha^{2}+1}{2 \sqrt{5} } = \frac{\alpha}{2}, \quad \frac{\alpha^{6}+1}{2 \sqrt{5} } = \alpha^3, \quad 
\frac{\alpha^{10}+1}{5 \sqrt{5} } = \alpha^5, \quad \frac{\alpha^{18}+1}{17\sqrt{5}} = 2 \alpha^9.$$ 
Using these dependencies we obtain
\begin{align*}
{\it\Lambda}_B & = (n_2+1) \log \alpha - (a_2+1) \log 2, &&
{\it\Lambda}_B = (n_2+3) \log \alpha - a_2 \log 2, \\
{\it\Lambda}_B &= (n_2+5) \log \alpha - a_2 \log 2 \quad \mbox{and}\;
&&{\it\Lambda}_B = (n_2+9) \log \alpha - (a_2-1) \log 2
\end{align*}
for $(k,r) = (0, 2), (0, 6), (2, 10), (4,18)$ respectively. Thus we get 
$$\left| \gamma - \frac{a_2+1}{n_2 + 1} \right| <  \frac{5}{2^{a_1 - a_3} (n_2 + 1)}, \quad \quad \;\;
\left| \gamma - \frac{a_2  }{n_2 + 3} \right| <  \frac{5}{2^{a_1 - a_3} (n_2 + 3)},$$
$$\left| \gamma - \frac{a_2}{n_2 + 5} \right| <  \frac{5}{2^{a_1 - a_3} (n_2 + 5)} \quad \mbox{and}\quad 
\left| \gamma - \frac{a_2 - 1}{n_2 + 9} \right| <  \frac{5}{2^{a_1 - a_3} (n_2 + 9)} $$
respectively. If $a_1 - a_3 \leq 211 $ the previous bound is still true. Now assume $a_1 - a_3 > 211 $. Then $2^{a_1 - a_3} > 4.2 \cdot 10^{63} > 10 (n_2 + 9) $, hence
$$ \frac{5}{2^{a_1 - a_3} (n_2 + 1)} < \frac{1}{2 (n_2 + 1)^2}, \quad  \quad \;\;
  \frac{5}{2^{a_1 - a_3} (n_2 + 3)} < \frac{1}{2 (n_2 + 3)^2} ,$$
$$ \frac{5}{2^{a_1 - a_3} (n_2 + 5)} < \frac{1}{2 (n_2 + 5)^2} \quad \mbox{and}\quad 
  \frac{5}{2^{a_1 - a_3} (n_2 + 9)} < \frac{1}{2 (n_2 + 9)^2} $$
respectively. By a criterion of Legendre each of $\frac{a_2+1}{n_2 + 1}$, $\frac{a_2}{n_2 + 3}$, $\frac{a_2}{n_2 + 5}$ and $\frac{a_2-1}{n_2 + 9}$ is a convergent to $\gamma$
and we may assume that $\frac{a_2+1}{n_2 + 1}$, $\frac{a_2}{n_2 + 3}$, $\frac{a_2}{n_2 + 5}$ and $\frac{a_2-1}{n_2 + 9}$ is of the form $\frac{p_j}{q_j}$ for some $j=0, 1, 2, \dotso, 124$.
Indeed, we may assume that $j\leq 124$ since $q_{125} > 4.2 \cdot 10^{62} > n_2 + 9$ but $q_{124} < 4.2 \cdot 10^{62}$. However it is well known (see e.g. \cite[page 47]{Baker:NT}) that 
$$  \frac{1}{(s_{j+1} + 2)q_j^2} < \left| \gamma - \frac{p_j}{q_j} \right|.$$
and since $\max \lbrace s_{j+1}: j=0, 1, 2, \dotso, 124 \rbrace = 134$, we have
$$\frac{1}{136q_j^2} < \frac{5}{2^{a_1-a_3} q_j}$$
and $q_j$ divides one of $\lbrace n_2+1, n_2+3, n_2+5, n_2+9 \rbrace$.
Thus the inequality
$$ 2^{a_1-a_3} < 5 \cdot 136(n_2+9) < 5 \cdot 136 \cdot 4.2 \cdot 10^{62}$$
yields $a_1-a_3 < 218$. Hence even in the case that $(k, r) \in \lbrace (0, 2), (0, 6), (2, 10), (4, 18)\rbrace$ we obtain the upper bound  $a_1-a_3  \leq 233$. \\ 

\noindent \textbf{Step 5:} \textit{We consider Case 2 and show that under the assumption that $n_1-n_2 \leq 315$ we have that $a_1-a_2 \leq  224$.}\\

In this step we consider inequality~\eqref{Case2} and assume that $a_1-a_2, a_1 - a_3 \geq 20$. Recall that
$${\it\Lambda}_2 =  n_2 \log \alpha - a_1 \log 2 + \log \left(\frac{  \alpha^{n_1 - n_2} + 1 }{\sqrt{5}} \right)$$
and inequality \eqref{Case2} yields that 
$ |{\it\Lambda}_2| < 4.9 \cdot 2^{-(a_1 - a_2)}$. 
Then we get  
\[ 0< \left| n_2 \cdot \frac{\log \alpha}{\log 2}  -a_1 +
\frac{\log ((\alpha^{n_1 - n_2} + 1)/\sqrt{5}) }{\log 2} \right|
< 8 \cdot 2^{-(a_1 - a_2)}.\]
We apply our algorithm with the same $\gamma$ and $M$, but we use $(A,B)=(8, 2)$ and 
\[ \mu=\mu_r = \frac{\log ((\alpha^{r} + 1)/\sqrt{5}) }{\log 2} ,\]
for each possible value of $n_1-n_2=r=0,1,\dots, 315$.
Similar as in Step 4 we obtain $a_1-a_2 \leq 224$, except in the problematic case that $r \in \lbrace 2, 6 \rbrace$.
However these two problematic cases can be treated in a similar way as the problematic cases in Step 4. That is we find a multiplicative relation
between $2$, $\alpha$ and $ \frac{\log ((\alpha^{r} + 1)/\sqrt{5}) }{\log 2}$ and reduce linear form ${\it\Lambda}_2$ to a linear form in two logarithms
and use the theory of continued fractions to obtain also in these problematic cases upper bounds for $a_1-a_2$. Thus in any case we obtain $a_1-a_2 \leq 224$.\\ 

\noindent \textbf{Step 6:} \textit{We continue to consider Case 2 and show that under the assumption that $n_1-n_2 \leq 315$ and $a_1-a_2 \leq 224$ we have that $a_1-a_3 \leq  233$.}\\

Now we have $n_1 - n_2 \leq 315$ and $a_1-a_2 \leq 224$ and we shall assume that $a_1 - a_3 \geq 20$ and attempt to reduce the huge upper bound 
for $a_1 - a_3$ with the use of inequality \eqref{CaseB}. This setting has already been considered in Case 1B, where we obtained
$$ a_1 - a_3 \leq 233.$$ 

Table \ref{TableCase1Reduce} summarizes our results obtained so far.

\begin{table}[ht] \label{Table3}
\caption{Summary of results} \label{TableCase1Reduce} 
\centering  
\begin{tabular}{c c c c c}  
\hline  
\textbf{Upper bound of $(\leq)$}  & \textbf{Case 1A} &\textbf{ Case 1B} & \textbf{Case 2 } & \textbf{Overall } \\ [0.5ex]  
\hline  
$a_1 - a_2$ 
& $218$ & $218$ & $224$ & $224$ \\ 
$a_1 - a_3$ 
& $225$  & $233$  & $233$ & $233$ \\
$n_1 - n_2$ 
& $334$  & $324$ & $315$ & $334$ \\
[1ex]  
\hline  
\end{tabular}
\end{table}

\noindent \textbf{Step 7:} \textit{Under the assumption that $n_1-n_2 \leq 334$, $a_1-a_2 \leq 224$ and $a_1-a_3 \leq  233$ we show that $n_1\leq 343$.}\\

For the last step we consider inequality \eqref{Case3}. Recall that
$${\it\Lambda}_3 = - n_1 \log \alpha + a_1 \log 2 + \log \left( \frac{\sqrt{5} (1+ 2^{a_2-a_1} + 2^{a_3-a_1} )}  { 1 + \alpha^{n_2-n_1} } \right) $$
and inequality \eqref{Case3} yields that
$ |{\it\Lambda}_3| < 2.02 \alpha^{-n_1} $. 
Then we get
\[ 0< \left| n_1 \cdot \frac{\log \alpha}{\log 2}  -a_1 +
\frac{\log \left(\left(1 + \alpha^{n_2-n_1}\right)/\left(\sqrt{5} \left(1+ 2^{a_2-a_1} + 2^{a_3-a_1} \right)\right)\right) }{\log 2} \right|
< 3 \alpha^{-n_1}.\]
We proceed as described in Remark \ref{RemarkReduction} with the same $\gamma$ and $M$ as in the previous steps, but we use $(A,B)=(3, \alpha)$ and  
\[ \mu=\mu_{k,l,r} = \frac{\log \left(\left(1 + \alpha^{-r} \right)/\left(\sqrt{5} \left(1+ 2^{-k} + 2^{-l} \right)\right) \right) }{\log 2} ,\]
 for each possible value of $a_1 - a_2 = k = 0, 1, \dots, 224$ , $a_1 - a_3 = l = 0, 1, \dots, 233$ 
(with respect to the obvious condition that $a_1 - a_2 \leq a_1 - a_3$) and $n_1-n_2=r=0,1,\dots, 334$.
Starting with $q_{125}$ we compute the upper bound for $n_1$ by the formula
$ n_1 < \frac{\log(3q/\varepsilon)}{\log \alpha} $
for the respective choices of $q$ such that $\varepsilon>0$. For all triples $(k,l,r)$ except
$$(k, l, r) \in \lbrace  (0, 1, 10), (0, 3, 18), (1, 1, 2),  (1, 1, 6), (1, 3, 14), (3, 3, 10), (5, 5, 18) \rbrace$$
the algorithm terminates and yields  
\begin{equation} \label{nbound}
n_1 \leq 343. 
\end{equation} 
The problematic cases can be treated in a similar way as in Step 4 and yield similarly small upper bounds for $n_1$.
In particular we obtain that $n_1 \leq 343$ in all cases. However this upper bound contradicts our assumption that $n_1 \geq 360$.
Therefore no further solutions to \eqref{eqt3A_1} exist and Theorem \ref{Main3A_1} is proved.

\subsection{Proof of Theorem \ref{Main3A_2}}

We reduce the upper bound for $m_1$ obtained in Proposition~\ref{prop:bound3A_2} by applying several times our algorithm described in Remark \ref{RemarkReduction}.
We do this in a similar manner as in the proof of Theorem \ref{Main3A_1}.\\

\noindent \textbf{Step 1:} \textit{We show that $t_1-t_2 \leq 218$ or $m_1-m_2 \leq 314$.}\\

First, we consider inequality \eqref{Case0_2} and deduce that 
\begin{equation*} 
\begin{split}
0< \left| m_1 \cdot \frac{\log \alpha}{\log 2 }   -t_1 + \frac{\log(1/\sqrt{5})}{\log 2 } \right| 
< & \max\left\{ 43 \cdot 2^{-(t_1-t_2)} , 43 \alpha^{-(m_1-m_2)} \right\}.
\end{split}
\end{equation*}
We apply Lemma \ref{Reduction} with the same $\gamma= \frac{\log \alpha}{\log 2 }$ as in the case of Theorem \ref{eqt3A_1}, but we use $ M = 4.2 \cdot 10^{62} $,  
$(A,B)=(43, 2)$ or $(43, \alpha)$ and $ \mu = \frac{\log(1/\sqrt{5})}{\log 2}$. 
We consider the $125$-th convergent $ \frac{p_{125}}{q_{125}}$  of $\gamma$ and obtain $\varepsilon > 0.24$ and therefore either
\[ t_1-t_2 \leq \frac{\log(43q/0.24)}{\log 2} \leq 218, \quad \mbox{or} \quad m_1-m_2 \leq \frac{\log(43q/0.24)}{\log \alpha} \leq 314. \]

Now, we distinguish between 
\begin{description}
 \item [Case 1] $m_1-m_2 \leq 314$ \quad and
 \item [Case 2] $t_1-t_2 \leq 218$.
\end{description}

\noindent \textbf{Step 2:} \textit{We consider Case 1 and show that under the assumption that $m_1-m_2 \leq 314$ we have $ t_1-t_2 \leq 226 $ or $m_1-m_3 \leq 326$.}\\

We consider inequality~\eqref{Case1_2} and get
\[ 0< \left| m_2 \cdot\frac{\log \alpha}{\log 2 }   -t_1 +
\frac{\log \left(\left(\alpha^{ m_1-m_2 } + 1 \right)/{\sqrt{5} } \right) }{\log 2} \right|
< 36 \max \left\{ 2^{-(t_1 - t_2)}, \alpha^{-(m_1-m_3)} \right\}.\]
We apply our algorithm (cf. Remark \ref{RemarkReduction}) for each possible value of $m_1-m_2=k \leq 314$ and the algorithm yields $ t_1-t_2 \leq 226 $ or $m_1-m_3 \leq 326$
for all $k= 1,2,\dots, 314$ except $k \in \lbrace 2, 6 \rbrace$. These two problematic cases can be treated by using continued fractions and Legendre's criterion. Thus we obtain
in all cases that $ t_1-t_2 \leq 226 $ or $m_1-m_3 \leq 326$.

Within Case 1, we distinguish between the following two sub-cases:
\begin{description}
 \item [Case 1A] $m_1-m_3 \leq 326$ \quad and
 \item [Case 1B] $t_1-t_2 \leq 226$.
\end{description}

\noindent \textbf{Step 3:} \textit{We consider Case 1A and show that under the assumption that $m_1-m_2 \leq 314$ and $m_1-m_3 \leq 326$ we have $ t_1-t_2 \leq 231 $.}\\

We consider inequality~\eqref{CaseA_2} and get   
\[ 0< \left| m_1 \cdot\frac{\log \alpha}{\log 2 } - t_1  + \frac{\log \left(\left( 1 + \alpha^{m_2-m_1} + \alpha^{m_3-m_1}  \right)  / \sqrt{5} \right) }{\log 2} \right|
< 6 \cdot 2^{-(t_1-t_2)}.\] 
For each possible value of $m_1-m_2=k \leq 314$ and $m_1-m_3=l \leq  326$ (with respect to the obvious condition $m_1 - m_2 \leq m_1 - m_3$) except for
 $$(k, l) \in \lbrace (0, 3), (1, 1), (1, 5), (3, 4), (7, 8) \rbrace, $$
our algorithm yields $t_1-t_2 \leq 231$. Note that the same upper bound can be concluded for the exceptional cases by using continued fractions and Legendre's criterion.\\

\noindent \textbf{Step 4:} \textit{We consider Case 1B and show that under the assumption that $m_1-m_2 \leq 314$ and $t_1-t_2 \leq 226$ we have $m_1-m_3 \leq 336 $.}\\

In view of Step 6 we consider $m_1 - m_2 \leq 323$ instead of $m_1 - m_2 \leq 314$ as required in this step. 
We consider inequality~\eqref{CaseB_2} and get 
\[ 0< \left| m_2 \cdot\frac{\log \alpha}{\log 2 }- t_2  + \frac{\log \left( (\alpha^{m_1 - m_2} + 1) / \left( \sqrt{5} \left(2^{t_1-t_2} + 1 \right) \right)\right) }{\log 2} \right|
< 9 \alpha^{-(m_1-m_3)} .\] 
By applying our algorithm for each possible value of $m_1-m_2=k \leq 323$ and $t_1-t_2=r  \leq 226$  we get
$  m_1-m_3  \leq 336$, except for $(k, r) \in \lbrace (2, 0), (6, 0), (10, 2), (18, 4)\rbrace$). However, by using continued fractions and Legendre's criterion
we obtain the same upper bound also for these exceptional cases. \\

\noindent \textbf{Step 5:} \textit{We consider Case 2 and show that under the assumption that $t_1-t_2 \leq 218$ we have $m_1-m_2 \leq 323$.}\\

We consider inequality~\eqref{Case2_2} and get 
\[ 0< \left| m_2 \cdot\frac{\log \alpha}{\log 2 } -t_2 +
\frac{\log \left(1/ \left(\sqrt{5} \left(2^{t_1 - t_2} + 1\right) \right)\right) }{\log 2} \right|
< 12 \alpha^{-(m_1 - m_2)}.\] 
For each possible value of $t_1-t_2=r  \leq  218$ our algorithm yields $ m_1-m_2 \leq 323$.\\

\noindent \textbf{Step 6:} \textit{We continue to consider Case 2 and show that under the assumption that $t_1-t_2 \leq 218$ and $m_1-m_2 \leq 323$ we have $ m_1 - m_3 \leq 336$.}\\

This situation is covered by Step 4 and we obtain that $ m_1 - m_3 \leq 336$. Table \ref{TableCase2Reduce} summarizes our results obtained so far.

\begin{table}[ht] 
\caption{Summary of results} \label{TableCase2Reduce} 
\centering  
\begin{tabular}{c c c c c}  
\hline 
\textbf{Upper bound of $(\leq)$}  & \textbf{Case 1A} &\textbf{ Case 1B} & \textbf{Case 2 } & \textbf{Overall }\\ [0.5ex] 
\hline 
$m_1 - m_2$ 
&  $314$ &  $314$ &  $323$  &  $323$\\ 
$m_1 - m_3$ 
& $326$  &  $336$  &  $336$ &  $336$\\
$t_1 - t_2$ 
& $231$   &  $226$  & $218$ & $231$ \\
[1ex] 
\hline 
\end{tabular}
\end{table}

\noindent \textbf{Step 7:} \textit{Under the assumption that $t_1-t_2 \leq 231$, $m_1-m_2 \leq 323$ and $m_1-m_3 \leq  336$ we show that $m_1\leq 353$.}\\

For the last step in our reduction process we consider inequality \eqref{Case3_2} and get  
\[ 0< \left| m_1\cdot\frac{\log \alpha}{\log 2 }  -t_1 +
\frac{\log \left( \left(1+ \alpha^{m_2-m_1} + \alpha^{m_3-m_1} \right)/\left(\sqrt{5} \left(1 + 2^{t_2-t_1}\right) \right)\right) }{\log 2} \right|
< 6 \alpha^{-m_1}.\]
We apply our algorithm for each possible value of $m_1 - m_2 = k \leq 324$ , $m_1 - m_3 = l \leq 337$ (with respect to the obvious condition $m_1 - m_2 \leq m_1 - m_3$) and $t_1-t_2=r \leq  232 $
and get $m_1 \leq 353$ except in the case that 
$$
(k, l, r) \in \lbrace (0, 3, 0), (1, 1, 0), (1, 5, 0), (3, 4, 0), (7, 8, 0), (1,9,2), (11,12,2), (1,17,4), (19,20,4)  \rbrace.
$$ 
These exceptional cases can be treated by using continued fractions and Legendre's criterion. Thus we obtain the upper bound  $m_1 \leq 353$ in all cases. 
But this upper bound contradicts our assumption that $m_1 \geq 360$. Therefore, no further solutions to \eqref{eqt3A_2} exist and Theorem \ref{Main3A_2} is proved. 

\section{Appendix - Lists of solutions for Theorem \ref{Main3A_1} and Theorem \ref{Main3A_2}}

The solutions for Diophantine equation \eqref{eqt3A_1} in Theorem \ref{Main3A_1} are displayed below. Since $F_1 = F_2$, the solutions involving $F_1$ are not displayed for the sake of simplicity. 
\begin{align*}
&  F_{3} + F_{2} = 2^{0} + 2^{0} + 2^{0} = 3,
&& F_{3} + F_{3} = 2^{1} + 2^{0} + 2^{0} = 4, \\
& F_{4} + F_{0} = 2^{0} + 2^{0} + 2^{0} = 3,
&& F_{4} + F_{2} = 2^{1} + 2^{0} + 2^{0} = 4, \\
&  F_{4} + F_{3} = 2^{1} + 2^{1} + 2^{0} = 5,
&& F_{4} + F_{4} = 2^{1} + 2^{1} + 2^{1} = 6, \\
&  F_{4} + F_{4} = 2^{2} + 2^{0} + 2^{0} = 6,
&& F_{5} + F_{0} = 2^{1} + 2^{1} + 2^{0} = 5,\\
& F_{5} + F_{2} = 2^{1} + 2^{1} + 2^{1} = 6, 
&& F_{5} + F_{2} = 2^{2} + 2^{0} + 2^{0} = 6, \\
&  F_{5} + F_{3} = 2^{2} + 2^{1} + 2^{0} = 7, 
&& F_{5} + F_{4} = 2^{2} + 2^{1} + 2^{1} = 8, \\
& F_{ 5 } + F_{ 5 } = 2^{ 2 } + 2^{ 2 } + 2^{ 1 } = 10 , 
&& F_{5} + F_{5} = 2^{3} + 2^{0} + 2^{0} = 10,  \\
&  F_{6} + F_{0} = 2^{2} + 2^{1} + 2^{1} = 8,
&& F_{6} + F_{2} = 2^{2} + 2^{2} + 2^{0} = 9,\\
& F_{ 6 } + F_{3  } = 2^{ 2 } + 2^{ 2 } + 2^{ 1 } = 10, 
&& F_{6} + F_{3} = 2^{3} + 2^{0} + 2^{0} = 10,\\
& F_{6} + F_{4} = 2^{3} + 2^{1} + 2^{0} = 11,
&&  F_{6} + F_{5} = 2^{3} + 2^{2} + 2^{0} = 13,  \\
& F_{ 6 } + F_{ 6 } = 2^{ 3 } + 2^{ 2 } + 2^{ 2 } = 16, 
&& F_{7} + F_{0} = 2^{3} + 2^{2} + 2^{0} = 13, \\
& F_{ 7 } + F_{2  } = 2^{ 3 } + 2^{ 2 } + 2^{ 1 } = 14, 
&& F_{ 7 } + F_{ 4 } = 2^{ 3 } + 2^{ 2 } + 2^{ 2 } = 16,  \\
& F_{ 7 } + F_{5  } = 2^{ 3 } + 2^{ 3 } + 2^{ 1 } = 18 , 
&& F_{7} + F_{5} = 2^{4} + 2^{0} + 2^{0} = 18, \\
& F_{7} + F_{6} = 2^{4} + 2^{2} + 2^{0} = 21,
&&  F_{ 7 } + F_{ 7 } = 2^{ 4 } + 2^{ 3 } + 2^{ 1 } = 26 ,  \\
& F_{8} + F_{0} = 2^{4} + 2^{2} + 2^{0} = 21,
&& F_{ 8 } + F_{2  } = 2^{ 4 } + 2^{ 2 } + 2^{ 1 } = 22 , \\
& F_{ 8 } + F_{ 4 } = 2^{ 3 } + 2^{ 3 } + 2^{3 } = 24 ,
&& F_{ 8 } + F_{4  } = 2^{ 4 } + 2^{ 2 } + 2^{ 2 } = 24 , \\
& F_{ 8 } + F_{ 5 } = 2^{ 4 } + 2^{ 3 } + 2^{1 } = 26 ,
&& F_{ 8 } + F_{ 7  } = 2^{ 4 } + 2^{ 4 } + 2^{ 1 } = 34 ,\\
&  F_{8} + F_{7} = 2^{5} + 2^{0} + 2^{0} = 34,
&& F_{ 8 } + F_{ 8 } = 2^{ 5 } + 2^{ 3 } + 2^{1 } = 42 , \\
& F_{9} + F_{0} = 2^{4} + 2^{4} + 2^{1} = 34,
&& F_{9} + F_{0} = 2^{5} + 2^{0} + 2^{0} = 34, \\
&  F_{9} + F_{2} = 2^{5} + 2^{1} + 2^{0} = 35,
&& F_{ 9 } + F_{ 3  } = 2^{ 4 } + 2^{ 4 } + 2^{ 2 } = 36 ,  \\
& F_{ 9 } + F_{ 3 } = 2^{ 5 } + 2^{ 1 } + 2^{1 } = 36 ,
&& F_{9} + F_{4} = 2^{5} + 2^{2} + 2^{0} = 37,  \\
& F_{ 9 } + F_{ 6  } = 2^{ 5 } + 2^{ 3 } + 2^{ 1 } = 42 , 
&& F_{ 9 } + F_{ 9 } = 2^{ 5 } + 2^{ 5 } + 2^{2 } = 68 , \\
& F_{ 9 } + F_{ 9  } = 2^{ 6 } + 2^{ 1 } + 2^{ 1 } = 68 , 
&& F_{ 10 } + F_{ 2 } = 2^{ 5 } + 2^{ 4 } + 2^{3 } = 56 , \\
& F_{ 10 } + F_{ 7  } = 2^{ 5 } + 2^{ 5 } + 2^{ 2 } = 68 ,
&& F_{ 10 } + F_{ 7 } = 2^{ 6 } + 2^{ 1 } + 2^{1 } = 68 , \\
& F_{ 10 } + F_{ 8  } = 2^{ 6 } + 2^{ 3 } + 2^{ 2 } = 76 ,
&& F_{11} + F_{6} = 2^{6} + 2^{5} + 2^{0} = 97, \\
&  F_{ 11 } + F_{ 10 } = 2^{ 6 } + 2^{ 6 } + 2^{4 } = 144 ,
&& F_{ 11 } + F_{ 10  } = 2^{ 7 } + 2^{ 3 } + 2^{ 3 } = 144 , \\
&  F_{12} + F_{0} = 2^{6} + 2^{6} + 2^{4} = 144,
&& F_{12} + F_{0} = 2^{7} + 2^{3} + 2^{3} = 144,\\
& F_{12} + F_{2} = 2^{7} + 2^{4} + 2^{0} = 145,
&& F_{ 12 } + F_{ 3} = 2^{ 7 } + 2^{ 4 } + 2^{1 } = 146 , \\
&  F_{ 12 } + F_{ 6  } = 2^{ 7 } + 2^{ 4 } + 2^{ 3 } = 152 ,
&& F_{ 12 } + F_{ 12 } = 2^{ 7 } + 2^{ 7 } + 2^{ 5 } = 288 , \\
& F_{ 12 } + F_{ 12  } = 2^{ 8 } + 2^{ 4 } + 2^{ 4 } = 288 , 
&& F_{ 13 } + F_{ 10 } = 2^{ 7 } + 2^{ 7 } + 2^{ 5 } = 288 , \\
&  F_{ 13 } + F_{ 10  } = 2^{ 8 } + 2^{ 4 } + 2^{ 4 } = 288 , 
&& F_{ 13 } + F_{ 11 } = 2^{ 8 } + 2^{ 6 } + 2^{ 1 } = 322 , \\
& F_{14} + F_{6} = 2^{8} + 2^{7} + 2^{0} = 385,
&& F_{14} + F_{12} = 2^{9} + 2^{3} + 2^{0} = 521, \\
& F_{ 15 } + F_{ 9  } = 2^{ 9 } + 2^{ 7 } + 2^{ 2 } = 644 , 
&& F_{ 16 } + F_{ 10 } = 2^{ 10 } + 2^{ 4 } + 2^{ 1 } = 1042 , \\
& F_{ 17 } + F_{ 4  } = 2^{ 10 } + 2^{ 9 } + 2^{ 6 } = 1600 ,
&& F_{ 18 } + F_{ 6 } = 2^{ 11 } + 2^{ 9 } + 2^{ 5 } = 2592.
\end{align*}

The solutions for Diophantine equation \eqref{eqt3A_2} in Theorem \ref{Main3A_2} are displayed below. Since $F_1 = F_2$, the solutions involving $F_1$ are not displayed for the sake of simplicity. 
\begin{align*}
& F_{  2 } + F_{   2  } + F_{ 0} = 2^{ 0} + 2^{0} =2, 
&& F_{ 2 } + F_{   2  } + F_{ 2} = 2^{ 1} + 2^{0} =3,  \\
& F_{ 3 } + F_{  0  } + F_{0 }  = 2^{0} + 2^{0} = 2,
&& F_{ 3 } + F_{   2  } + F_{ 0} = 2^{ 1} + 2^{0} =3, \\
& F_{ 3 } + F_{ 2  } + F_{ 2  } = 2^{ 1 } + 2^{ 1 } = 4 ,
&& F_{ 3 } + F_{   3  } + F_{ 0} = 2^{ 1} + 2^{1 } =4,\\
&  F_{ 3 } + F_{  3   } + F_{ 2} = 2^{ 2} + 2^{ 0 } = 5, 
&& F_{ 3} + F_{ 3 } + F_{ 3  }  = 2^{ 2 } + 2^{ 1 } = 6 , \\
&  F_{ 4 } + F_{ 0  } + F_{0 }  = 2^{1} + 2^{0} = 3, 
&& F_{ 4 } + F_{   2  } + F_{ 0} = 2^{ 1} + 2^{1 } = 4,  \\
&  F_{ 4 } + F_{   2  } + F_{2} = 2^{ 2} + 2^{0 } = 5, 
&& F_{ 4 } + F_{ 3  } + F_{0 }  = 2^{2} + 2^{0} = 5,\\
&  F_{ 4 } + F_{ 3  } + F_{ 2  } = 2^{ 2 } + 2^{ 1 } = 6 ,  
&& F_{ 4 } + F_{ 4} + F_{0 }  = 2^{2} + 2^{1} = 6,  \\
&  F_{ 4 } + F_{ 4 } + F_{ 3  }  = 2^{ 2 } + 2^{ 2 } = 8 , 
&& F_{ 4 } + F_{ 4  } + F_{ 4 } = 2^{ 3} + 2^{0 } = 9,  \\
& F_{ 5 } + F_{ 0} + F_{0 }  = 2^{2} + 2^{0} = 5,
&& F_{ 5 } + F_{ 2} + F_{0 }  = 2^{2} + 2^{1} = 6,\\
&  F_{ 5 } + F_{ 3  } + F_{ 2  } = 2^{ 2 } + 2^{ 2 } = 8 , 
&& F_{ 5 } + F_{ 3 } + F_{3  }  = 2^{3} + 2^{0} = 9,\\
&   F_{ 5 } + F_{  4   } + F_{ 0  } = 2^{ 2} + 2^{2 } = 8, 
&& F_{ 5 } + F_{  4  } + F_{ 2  } = 2^{ 3} + 2^{0 } = 9, \\
&   F_{ 5 } + F_{ 4 } + F_{ 3  }  = 2^{ 3 } + 2^{ 1 } = 10 ,
&& F_{ 5 } + F_{ 5  } + F_{ 0  }  = 2^{3} + 2^{1} = 10,  \\
&  F_{ 5 } + F_{ 5  } + F_{ 3  } = 2^{ 3 } + 2^{ 2 } = 12 , 
&& F_{ 6 } + F_{ 0  } + F_{ 0  } = 2^{ 2} + 2^{2 } = 8,\\
& F_{ 6} + F_{ 2  } + F_{ 0  }  = 2^{3} + 2^{0} = 9,
&& F_{ 6 } + F_{ 2 } + F_{ 2  }  = 2^{ 3 } + 2^{ 1 } = 10, \\
&  F_{ 6} + F_{  3  } + F_{ 0  }  = 2^{3} + 2^{ 1} = 10,  
&& F_{ 6 } + F_{ 3  } + F_{ 3  } = 2^{ 3 } + 2^{ 2 } = 12 ,\\
&  F_{ 6 } + F_{ 4 } + F_{ 2  }  = 2^{ 3 } + 2^{ 2 } = 12 , 
&&  F_{ 6 } + F_{ 5  } + F_{ 4  } = 2^{ 3 } + 2^{ 3 } = 16 ,  \\
& F_{ 6 } + F_{ 5 } + F_{ 5  }  = 2^{ 4 } + 2^{ 1 } = 18 ,
&&  F_{ 6} + F_{  6  } + F_{ 0  }  = 2^{3} + 2^{ 3} = 16,\\
&  F_{ 6} + F_{  6  } + F_{ 2  }  = 2^{4} + 2^{ 0} = 17, 
&& F_{ 6 } + F_{ 6  } + F_{ 3  } = 2^{ 4 } + 2^{ 1 } = 18 , \\
&  F_{ 6 } + F_{ 6 } + F_{ 6  }  = 2^{ 4 } + 2^{ 3 } = 24 ,
&& F_{ 7 } + F_{ 3  } + F_{ 2  } = 2^{ 3 } + 2^{ 3 } = 16 , \\
&  F_{ 7 } + F_{ 3  } + F_{ 3  }  = 2^{4} + 2^{ 0} = 17,
&& F_{ 7 } + F_{   4  } + F_{ 0  } = 2^{ 3} + 2^{3 } = 16,  \\
& F_{ 7 } + F_{  4  } + F_{ 2  } = 2^{ 4 } + 2^{0 } = 17,
&& F_{ 7 } + F_{ 4 } + F_{ 3  }  = 2^{ 4 } + 2^{ 1 } = 18 , \\
&  F_{ 7 } + F_{  5  } + F_{ 0  }  = 2^{4} + 2^{ 1 } = 18,
&& F_{ 7 } + F_{ 5  } + F_{ 3  } = 2^{ 4 } + 2^{ 2 } = 20 , \\
&   F_{ 7 } + F_{ 6 } + F_{ 4  }  = 2^{ 4 } + 2^{ 3 } = 24 ,
&& F_{ 7 } + F_{ 7  } + F_{ 6  } = 2^{ 5 } + 2^{ 1 } = 34 , \\
&  F_{ 8 } + F_{ 3 } + F_{ 2  }  = 2^{ 4 } + 2^{ 3 } = 24 ,
&& F_{ 8 } + F_{  4 } + F_{ 0  }  = 2^{4} + 2^{ 3 } = 24, \\
& F_{ 8 } + F_{ 6 } + F_{ 4  } = 2^{ 4 } + 2^{ 4 } = 32 ,
&& F_{ 8 } + F_{ 6 } + F_{ 5  }  = 2^{ 5 } + 2^{ 1 } = 34 , \\
&  F_{ 8 } + F_{ 7  } + F_{ 0  } = 2^{ 5 } + 2^{ 1 } = 34,
&& F_{ 8 } + F_{ 7 } + F_{ 3 } = 2^{ 5 } + 2^{ 2 } = 36 , \\
&   F_{ 9} + F_{ 0} + F_{ 0  }  = 2^{5} + 2^{ 1 } = 34,
&& F_{ 9 } + F_{ 2 } + F_{ 2  }  = 2^{ 5 } + 2^{ 2 } = 36 , \\
&  F_{ 9} + F_{ 3} + F_{ 0  }  = 2^{5} + 2^{ 2 } = 36,
&& F_{ 9 } + F_{ 4 } + F_{ 4 } = 2^{ 5 } + 2^{ 3 } = 40 ,\\
&  F_{ 9 } + F_{ 5 } + F_{ 2  }  = 2^{ 5 } + 2^{ 3 } = 40 ,
&& F_{ 9 } + F_{ 7 } + F_{ 2 } = 2^{ 5 } + 2^{ 4 } = 48 ,  \\
& F_{ 9 } + F_{ 8 } + F_{ 7  }  = 2^{ 6 } + 2^{ 2 } = 68 , 
&& F_{ 9 } + F_{  9} + F_{ 0  } = 2^{ 6 } + 2^{ 2 } = 68,\\
&  F_{ 10} + F_{ 5 } + F_{ 5  }  = 2^{6} + 2^{ 0 } = 65,
&& F_{ 10 } + F_{ 6 } + F_{ 2 } = 2^{ 5 } + 2^{ 5 } = 64 , \\
& F_{ 10} + F_{ 6 } + F_{ 3  }  = 2^{6} + 2^{ 0 } = 65, 
&& F_{ 10 } + F_{ 6 } + F_{ 4  }  = 2^{ 6 } + 2^{ 1 } = 66 , \\
& F_{ 10 } + F_{ 6 } + F_{ 5 } = 2^{ 6 } + 2^{ 2 } = 68 ,
&& F_{ 10 } + F_{ 7  } + F_{ 0 } = 2^{ 6 } + 2^{ 2 } = 68, \\
& F_{ 10 } + F_{ 10 } + F_{ 9  }  = 2^{ 7 } + 2^{ 4 } = 144 , 
&& F_{ 11 } + F_{ 5 } + F_{ 3 } = 2^{ 6 } + 2^{ 5 } = 96 ,\\
&  F_{ 11 } + F_{ 9 } + F_{ 5  }  = 2^{ 6 } + 2^{ 6 } = 128 ,
&& F_{ 11 } + F_{ 9 } + F_{ 7 } = 2^{ 7 } + 2^{ 3 } = 136 ,  \\
& F_{ 11 } + F_{ 9 } + F_{ 8  }  = 2^{ 7 } + 2^{ 4 } = 144 ,
&&  F_{ 11} + F_{ 10 } + F_{ 0  }  = 2^{ 7} + 2^{ 4 } = 144,  \\
& F_{ 12 } + F_{ 0  } + F_{ 0 } = 2^{ 7 } + 2^{ 4 } = 144 ,
&& F_{ 12 } + F_{ 6 } + F_{ 6 } = 2^{ 7 } + 2^{ 5 } = 160 , \\
& F_{ 12 } + F_{ 7 } + F_{ 4  }  = 2^{ 7 } + 2^{ 5 } = 160 ,
&&  F_{ 12 } + F_{ 11 } + F_{ 10 } = 2^{ 8 } + 2^{ 5 } = 288 , \\
& F_{ 12} + F_{ 12 } + F_{ 0  }  = 2^{ 8} + 2^{ 5 } = 288,
&&  F_{ 13 } + F_{ 8 } + F_{ 3  }  = 2^{ 7 } + 2^{ 7 } = 256 ,  \\
& F_{ 13 } + F_{ 8  } + F_{ 4  } = 2^{ 8 } + 2^{ 0 } = 257 ,
&&  F_{ 13 } + F_{ 9 } + F_{ 5 } = 2^{ 8 } + 2^{ 4 } = 272 , \\
&  F_{ 13 } + F_{ 9 } + F_{ 8  }  = 2^{ 8 } + 2^{ 5 } = 288 ,
&& F_{ 13} + F_{ 10 } + F_{ 0  }  = 2^{ 8} + 2^{ 5 } = 288,  \\
& F_{ 14 } + F_{ 5 } + F_{ 3 } = 2^{ 8 } + 2^{ 7 } = 384 , 
&& F_{ 14 } + F_{ 12 } + F_{ 10  }  = 2^{ 9 } + 2^{ 6 } = 576 ,\\
&  F_{ 16 } + F_{ 9 } + F_{ 4 } = 2^{ 9 } + 2^{ 9 } = 1024 , 
&& F_{ 16 } + F_{ 9 } + F_{ 5  }  = 2^{ 10 } + 2^{ 1 } = 1026 , \\
& F_{ 16 } + F_{ 12 } + F_{ 8 } = 2^{ 10 } + 2^{ 7 } = 1152 .\\
\end{align*}


\def\cprime{$'$}

\end{document}